\documentclass[11pt]{article}
\usepackage{amssymb,amsmath,amsthm,amsfonts}
\usepackage{bm}

\usepackage{graphicx}
\usepackage{subfig}
\usepackage{booktabs}
\usepackage{multirow}
\usepackage{array}
\usepackage{url}

\usepackage{cite}
\usepackage[colorlinks, linktocpage, citecolor=black, linkcolor=black]{hyperref}

\numberwithin{equation}{section}

\graphicspath{{figs/}{figs/designs/}{assets/}}

\begin{document}

% paper title
\title{\bf Agent4FreeFem: A Verification-Guided\\ Framework for Autonomous
Topology Optimization}
\author{Zhipu Cui\thanks{School of Mathematical Sciences, East China Normal
University, Shanghai 200241, China. (\texttt{52300155021@stu.ecnu.edu.cn})}
\and Helin Gong\thanks{SJTU Paris Elite Institute of Technology, Shanghai
Jiao Tong University, Shanghai 200240, China.
(\texttt{gonghelin@sjtu.edu.cn})}
\and Shengfeng Zhu\thanks{Corresponding author. School of Mathematical
Sciences, East China Normal University, Shanghai 200241, China.
(\texttt{sfzhu@math.ecnu.edu.cn})}}
\date{\today}

\maketitle

\begin{abstract}
Topology optimization has become an indispensable design tool
in the aerospace, automotive, and mechanical engineering industries. Yet bringing
a new problem from its description to an optimized result remains
slow, expert-dependent, and unreliable: the workflow spans many
manually connected stages, each new problem requires solver code to be
written and debugged manually, and such code often fails or yields
incorrect results. We present Agent4FreeFem, a pipeline driven by a large language model that turns a natural-language problem description into a set of verified optimization
results without manual solver-code development or debugging. The pipeline keeps the problem
definition unmodified, localizes errors using independent sources of
evidence, and completes constrained repairs automatically, thereby
providing efficient end-to-end automation, autonomous execution, and
verifiable results.
We evaluated Agent4FreeFem on a suite of single-physics and coupled
multiphysics problems. Every case passed runtime verification and
finished within the repair budget of each stage.
\par\vspace{1ex}

\vspace{0.25cm}
\end{abstract}

% =====================================================================
\section{Introduction}

Topology optimization determines the optimal material layout inside a
prescribed design domain under given physical constraints, without
assuming any preconceived structure. Since the seminal works \cite{bendsoe1988cmame,zhou1991coc}, the method has grown
into an indispensable industrial design
tool~\cite{bendsoe2003topology}: it is used to lighten
aerospace brackets and aircraft components, to stiffen automotive
chassis and engine mounts, to shape mechanical parts, and to design
civil structures, among others. By routinely
saving material and reducing weight while improving structural
performance, topology optimization shortens design cycles and
lowers
manufacturing cost, and with the spread of additive manufacturing its
intricate optimized geometries have become increasingly manufacturable.

Despite this maturity, taking a \emph{new} problem from a description
to an optimized result remains slow, expert-dependent, and fragile.
First, the process is time-consuming. The workflow spans many
disconnected stages, including mathematical formulation, mesh generation,
weak-form derivation, solver configuration, and tuning of the
optimization loop, each handled in a separate tool;
repeated cross-stage handoffs account for a substantial portion of the turnaround time. Second, the
process requires a high level of expertise. Every new
problem requires a specialist who understands both the underlying
physics and numerical computing to hand-write and debug the solver
code. Third, the process suffers from low reliability.
Hand-written solver code often compiles yet fails at runtime or
converges to incorrect solutions because of errors in weak forms,
signs, units, boundary conditions, or parameter scales, and locating
such errors is itself a major effort. Recent benchmarks of LLM-generated
finite-element code report similarly low pass rates~\cite{fembench2025,feabench2025},
confirming that merely writing code is far from producing a correct,
runnable solver.

Large language models (LLMs) are capable of generating
code~\cite{chen2021codex}, yet scientific coding remains
challenging for them~\cite{tian2024scicode}, and their finite-element
code still exhibits low pass rates~\cite{fembench2025,feabench2025}.
For mechanics and engineering design, multi-agent LLM collaborations
have been shown to solve
mechanics problems~\cite{ni2024mechagents}, agentic LLMs have been
fine-tuned for
finite-element methods~\cite{allfem2026}, and LLM agents generate
CAD models from text~\cite{texttocad2025}. On the industrial side,
LLM-based virtual assistants are being embedded into commercial CAE
products such as AnsysGPT~\cite{ansysgpt2024}, and autonomous agent
systems built on reasoning-and-acting frameworks handle large-scale
industrial scheduling~\cite{seo2026scheduling}, highlighting the
practical demand for autonomous engineering pipelines.

Topology optimization has begun to attract LLM-based automation.
Liang
et al.~\cite{lmto2025} integrate large language models into the
conceptual design stage of topology optimization: the model interprets
the natural-language design requirements and supports the development
of an initial conceptual design, yet its assistance stops at the
pre-processing stage, since the downstream numerical pipeline, from the PDE
formulation and discretization to the solver implementation, is still
assembled by experts, and the model's output is not automatically
linked to an executed optimizer. Yang et al.~\cite{yang2026simpllm}
use an LLM as an online optimization controller that monitors the
convergence of a SIMP solver and adaptively selects its continuation
parameters, but the LLM operates \emph{inside} the existing solver that
has already been written and debugged by experts, so the dominant cost
of the workflow, writing and debugging the solver code
itself, remains unaddressed. Heyrani Nobari et al.~\cite{oat2025} train
a deep generative foundation model that directly predicts optimized
topologies in a shape- and resolution-free manner, which shifts the
cost to data: it requires a large corpus of solved problems for
training, inherits the distribution of that corpus, and produces
designs without a PDE-constrained solve to verify the resulting designs.
Lin et al.~\cite{lin2026tomaster} present TO-Master, a
tool-orchestrated LLM agent that assembles density-based
topology-optimization models from natural-language instructions and
drives a fixed JAX-FEM optimization core; the supported formulations
are pre-implemented in the core, and the mesh and the boundary
conditions are confirmed by the human user. Park and
Chung~\cite{park2026topoptagents} introduce a multi-agent framework,
TopOptAgents, demonstrated on benchmark problems, in
which a specification validator, a code team, and a design critic
refine the formulation, the FEniCS implementation, and the converged
designs through four interleaved loops, with the critic judging rendered
designs against a rubric of validity, consistency, convergence, and
design quality.

To address these challenges, we present Agent4FreeFem, an agentic
pipeline that turns a natural-language problem description into a set
of verified optimization results without manual solver-code
development or debugging. The pipeline consists of
four stages. During \emph{Initialization}, an example-selection step
first scores the problem files against a reference library and selects
the most relevant reference examples; the LLM then parses the
description into a structured problem specification and a
mesh specification, annotating each value to indicate whether it was
\emph{extracted} from the source or \emph{designed} to fill a
gap, using the selected reference specifications as contextual
guidance. \emph{Generation} has the LLM write the solver code in the
AutoFreeFEM domain-specific language~\cite{allaire2024autofreefem}
in accordance with that specification, aided by the selected reference
programs and meshes. \emph{Execution} runs the generated code through the
AutoFreeFEM library and the FreeFem++ solver~\cite{freefemmanual}
with an automatic parameter sweep, producing optimized designs.
\emph{Reflection} has two responsibilities: when an execution
produces no design snapshot at all, it classifies the failure as a
program error, repairs the code, and re-executes; when the sweep runs
to completion, it first performs an overall check of the designs and
then post-processes the results, for example by composing the design
screenshots and cleaning up the temporary files.

The central methodological principle of Agent4FreeFem is to preserve
the problem definition, localize errors using independent sources of
evidence, and perform constrained repairs automatically. This principle is implemented in three ways. First, every value in the specification carries
its extracted/designed source annotation, so extracted values stay
verbatim, designed values carry rationales, and the problem definition
is never altered by any later stage. Second, the deterministic checks,
the semantic review, and the real execution evidence locate errors
independently of one another. Third, each repair modifies only the
flagged entries, and cases that fail to pass within the repair budget
are set aside for human review rather than silently accepted. This approach offers
three main advantages:

\begin{itemize}
\item \textbf{Efficient connectivity.} Problem parsing, mesh
handling, code generation, solving, and parameter tuning are
automatically chained inside a single pipeline, from the
natural-language description all the way to the finished sweep;
no manual handoff between stages is needed, so the entire flow
proceeds efficiently.
\item \textbf{Process autonomy.} The modeling specification, the
selection of reference examples, and the code generation are all performed
autonomously by the LLM, and the pipeline automatically runs the
solver and performs the parameter sweeps, so that no manual solver-code
development is required. Each decision is recorded in an
operation log, and the automatic parameter sweeps provide supporting
evidence, making the results reproducible and auditable.
\item \textbf{Process verifiability.} The generated code undergoes
deterministic static checks and runtime verification. If a check fails,
the pipeline automatically classifies the error, repairs the code, and
reruns it until the run passes or a given repair budget is exhausted.
\end{itemize}

We validate Agent4FreeFem with the level-set method on 34 different
examples spanning 6 problem types: linear elasticity, nonlinear
elasticity, linear thermal conduction, nonlinear thermal conduction,
Stokes flow, and coupled thermo-elasticity. All 2{,}050 solver runs passed runtime verification:
2{,}048 runs from the 32 suite cases, 64 parameter combinations per
case, and two runs from the extension cases.

The remainder of this paper is organized as follows. Section~2
reviews the preliminaries. Section~3
presents the methodology of Agent4FreeFem. Section~4 reports the
experimental results. Section~5 concludes the paper.

% =====================================================================
\section{Preliminaries}

This section reviews the methodological and computational foundations
of this work:
the level-set method of topology optimization and the
FreeFem++/AutoFreeFEM software stack.

%\subsection{Topology Optimization by the Level-Set Method}

Let $D\subset\mathbb{R}^d\ (d=2,3)$ be a fixed, bounded, open design domain with a boundary $\partial D$. Let $\Omega\subset D$ be a subdomain to be optimized. The level-set method for topology
optimization~\cite{osher1988jcp,allaire2004levelset,wang2003levelset}
represents the design implicitly through
$\phi:D\rightarrow\mathbb{R}$
\begin{equation}
\begin{cases}
\phi < 0, & \text{in } \Omega,\\
\phi > 0, & \text{in } D\setminus\overline{\Omega},\\
\phi = 0, & \text{on } \Gamma,
\end{cases}
\end{equation}
where $D=\Omega\cup\Gamma\cup(D\setminus\overline{\Omega})$ with $\Gamma$ being the interface to be optimized between two subdomains.
The material indicator is defined directly in terms of $\phi$ as follows
\begin{equation}
\chi(\phi)=
\begin{cases}
1, & \phi<0,\\
\chi_0, & \phi>0,
\end{cases}
\end{equation}
where $0 < \chi_0 \ll 1$.

A topology optimization model based on the level-set method can be expressed in the following general form
\begin{equation}\label{lsmTOP}
\begin{gathered}
\min_{\phi}\ J(\phi)=\mathcal{J}(\phi)+\ell\mathcal{V}(\phi)\\
%+\ell\int_D \chi(\phi)\,\mathrm{d}x\\
%\int_D \chi(\phi)\,j\bigl(\mathcal{U}(\phi)\bigr)\,\mathrm{d}x
%+\ell\int_D \chi(\phi)\,\mathrm{d}x\\
\begin{array}{cl}
\text{s.t.} & \begin{cases}
\mathcal{E}\bigl(\mathcal{U};\phi\bigr)=0 & \text{in}\ D,\\
\mathcal{B}\bigl(\mathcal{U}\bigr)=0 & \text{on}\ \partial D,
\end{cases}
\end{array}
\end{gathered}
\end{equation}
where $J$ denotes the total objective, $\mathcal{U}$ is the scalar or vector state field, $\ell$ is the material-cost weight,
$\mathcal{E}$ is the state operator, and $\mathcal{B}$ is the boundary-condition operator. The volume functional is defined by
$$\mathcal{V}(\phi)=\int_D \chi(\phi)\,\mathrm{d}x.$$
Here, $\mathcal{J}(\phi)$ denotes a general performance functional determined by the state fields and the design. The
level-set method has
been applied in many fields related to geometric optimization and reconstruction (see, e.g., \cite{zhu2018levelset,zheng2024rbffd,challis2009stokes,hu2025inverse,tan2024dgtopology,fan2025moving}).

The level-set function is transported so that the material interface moves \cite{allaire2013meshevolution,allaire2014meshevolution}, and the
finite-element mesh is then rebuilt to conform to the
updated interface, providing more accurate interface resolution than
traditional fixed-mesh level-set methods \cite{wang2003levelset,allaire2004levelset}. At each optimization iteration, an $H^1$-gradient flow is introduced by seeking a regularized descent field
$\boldsymbol{\theta}\in {V}:=[H^1(D)]^d$ such that
\begin{equation}
\int_D \bigl[
\boldsymbol{\theta}\cdot\mathbf{w}
+\alpha(\nabla\cdot\boldsymbol{\theta})
       (\nabla\cdot\mathbf{w})
+2\alpha\,
\boldsymbol{\varepsilon}(\boldsymbol{\theta})
:\boldsymbol{\varepsilon}(\mathbf{w})
\bigr]\,\mathrm{d}x
=-J'(\Omega)(\mathbf{w})
\quad \forall\,\mathbf{w}\in{V},
\end{equation}
where
\begin{equation}
\boldsymbol{\varepsilon}(\mathbf{w})
=\frac{1}{2}
\bigl[\nabla\mathbf{w}
+(\nabla\mathbf{w})^{T}\bigr]
\end{equation}
is the symmetric gradient, $\alpha>0$ is the diffusion
coefficient, and $J'(\Omega)(\mathbf{w})$ denotes the Eulerian derivative of the objective at the current $\Omega$ determined by $\phi$ in the direction $\mathbf{w}$ \cite{sokolowski1992shape}. Boundary constraints
are additionally enforced through prescribed
displacements and penalty terms on normal boundary
motion.

To seek a local minimizer of \eqref{lsmTOP}, we evolve a time-dependent level-set function $\phi=\phi(t,\mathbf{x}):\mathbb{R}^+\times D\rightarrow\mathbb{R}$ according to
\begin{equation}
\frac{\partial\phi}{\partial t}
+\boldsymbol{\theta}\cdot\nabla\phi=0
\quad \text{in}\ \mathbb{R}^+\times D,
\end{equation}
subject to an initial condition and homogeneous Neumann boundary conditions, where $t>0$ denotes the evolution time.
The transport equation is handled by a
characteristic-based advection solver \cite{yaji2016convected}. A first-order
approximation of the characteristic update from step $k\in\mathbb{N}$ to $k+1$ is
\begin{equation}
\widetilde{\phi}^{\,k+1}(\mathbf{x})
\approx
\phi^k\bigl(
\mathbf{x}-\tau_k\boldsymbol{\theta}^k(\mathbf{x})
\bigr),
\end{equation}
where the step size $\tau_k>0$ controls the magnitude of
the interface displacement.
Subsequently, the finite-element mesh is rebuilt to
follow the updated interface: the state problem is discretized on the updated geometry, and the
remesher~\cite{dapogny2014mmg} keeps the element edges conforming to the
interface and leaves the level-set function unchanged. A
signed-distance function is then reconstructed on
the new mesh using the fast marching method.
The reconstructed function approximately satisfies
the Eikonal condition
\begin{equation}
\begin{cases}
\phi^{k+1}=0,
& \text{on }\Gamma^{k+1},\\
|\nabla\phi^{k+1}|\approx 1,
& \text{away from }\Gamma^{k+1},
\end{cases}
\end{equation}
where $\Gamma^{k+1}$ denotes the updated material
interface. This reconstruction restores the
signed-distance representation without separately
integrating a time-dependent reinitialization equation.

The resulting level-set update can be summarized as follows
\begin{equation}
\begin{cases}
\widetilde{\phi}^{\,k+1}
=\mathcal{A}_{\tau_k\boldsymbol{\theta}^k}(\phi^k),\\
\mathcal{T}^{k+1}
=\mathcal{M}\bigl(\mathcal{T}^k,\widetilde{\phi}^{\,k+1}\bigr),\\
\phi^{k+1}
=\mathcal{R}\bigl(\mathcal{T}^{k+1}\bigr),
\end{cases}
\end{equation}
where $\mathcal{A}_{\tau_k\boldsymbol{\theta}^k}$ denotes level-set advection,
$\mathcal{M}$ denotes the remeshing
operator, $\mathcal{T}^k$ denotes the mesh at
iteration $k$, and $\mathcal{R}$ denotes
signed-distance reconstruction. The physical state
is subsequently
recomputed on the updated finite-element mesh, and
the objective value is evaluated. A backtracking
line-search procedure reduces the trial advection
step when the objective fails to decrease. This
process is repeated throughout the optimization.
Since the objective is re-evaluated on the remeshed
domain, the recorded objective value may increase. The
acceptance criterion requires that the new value be
lower than that of the previous iteration; if no tentative
step within the line-search budget decreases it, the
smallest tried step is taken. The optimization
always runs the full iteration budget, and the
reported design is the one attaining the smallest
objective value over the whole history of the run.

%\subsection{FreeFem++ and AutoFreeFEM}

The widely used finite-element environment FreeFem++~\cite{freefemmanual} allows a user to describe a
PDE in a high-level language close to its weak formulation, while
discretization, assembly, and solving are handled by the library.

AutoFreeFEM~\cite{allaire2024autofreefem}
raises the description level further: it is a symbolic Python layer on
top of FreeFem++ in which the user symbolically defines the domain, fields,
constants, objective functional, and weak forms, and the
library generates the FreeFem++ optimization scripts together with
\LaTeX{} documentation of the problem. Its central object is the
Lagrangian, assembled from the state fields, the test fields, the
objective $J$, and one weak-form expression. A companion
optimization-setup routine then configures the level-set optimizer with
the mesh, the mesh-size bounds, the 
gradient-flow diffusion coefficient and the advection step settings, the number of design iterations, the initial
level-set function, and the boundary-label lists.

% =====================================================================
\section{Methodology}

\begin{figure}[!tb]
\centering
\includegraphics[width=0.96\textwidth]{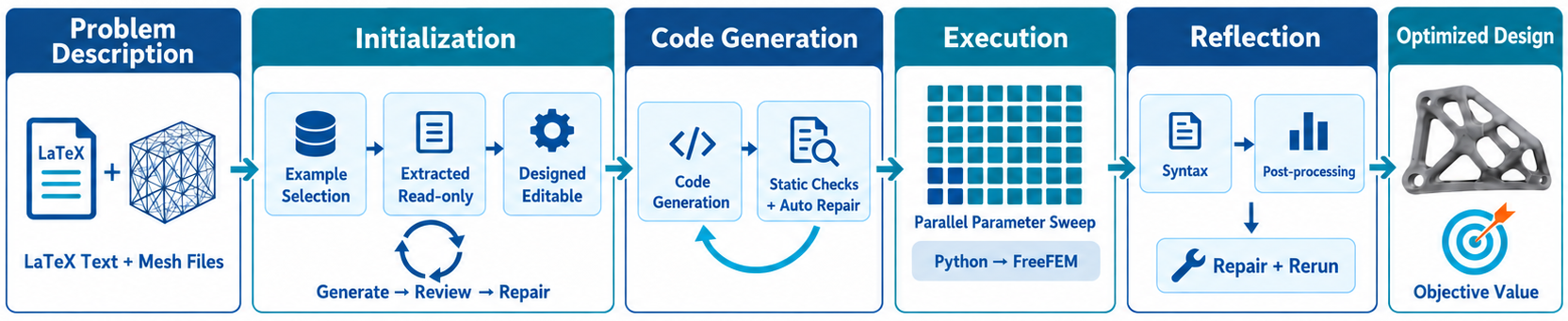}
\caption{Overall flow of Agent4FreeFem: from a natural-language
description through the four stages to the optimized designs.}
\label{fig:general}
\end{figure}

Agent4FreeFem takes as input a natural-language description of a
topology optimization problem and produces candidate optimized designs
for automated verification and subsequent selection.
Figure~\ref{fig:general} gives an overview of the pipeline. Four stages run in one chain:
\textbf{Initialization} selects the reference examples and turns the
description into a formal specification,
\textbf{Generation} turns the specification into solver code, \textbf{Execution} solves
the problem, and \textbf{Reflection} acts as the last line of defense: an execution that produces no design snapshot is classified as a program error and repaired, and a completed sweep is overall-checked and post-processed. Two features cut across all four
stages. First, every stage consumes exactly the artifacts the previous
stage produced, so no intermediate result is ever transcribed or passed
on by hand, and the pipeline can proceed from description to the
optimized designs in one uninterrupted run. Second, every stage leaves a complete,
append-only record of its inputs, outputs, and decisions, and every
value in the formal specification carries a provenance annotation
stating whether it was taken from the source or introduced by the
model to fill unspecified details.

\subsection{Initialization}

The first stage turns the unstructured problem description into
a unified form. It begins with example selection: an example
selector scores all of the problem's input files against a reference
library and packs the specifications of the highest-scoring examples
above a prescribed threshold into the parser's context under a token
budget, so the parser can base its design choices on the patterns of
similar problems. The description
may be a short statement of
intent, a detailed technical document, or something in between, and may
come with a mesh script that fixes the geometry. The model reads the
description and, using the selected reference specifications as
contextual guidance, produces two structured specifications: a
problem specification covering the physics, the fields and their
approximation spaces, the material constants with their units, the
boundary conditions, the objective and the weak forms, and the
optimization settings; and a mesh specification covering the geometry,
the boundary parametrizations with their labels, and the mesh
resolution. Every entry in both specifications carries one of two
annotations: \emph{extracted}, when it is transcribed from the
source, or \emph{designed}, when the model completes what
the source leaves open, together with a short rationale. This dual
annotation is the backbone of the pipeline: extracted values are
immutable ground truth that no later stage may change, while designed
values remain explicitly marked as choices.

\begin{figure}[!tb]
\centering
\includegraphics[width=0.98\textwidth]{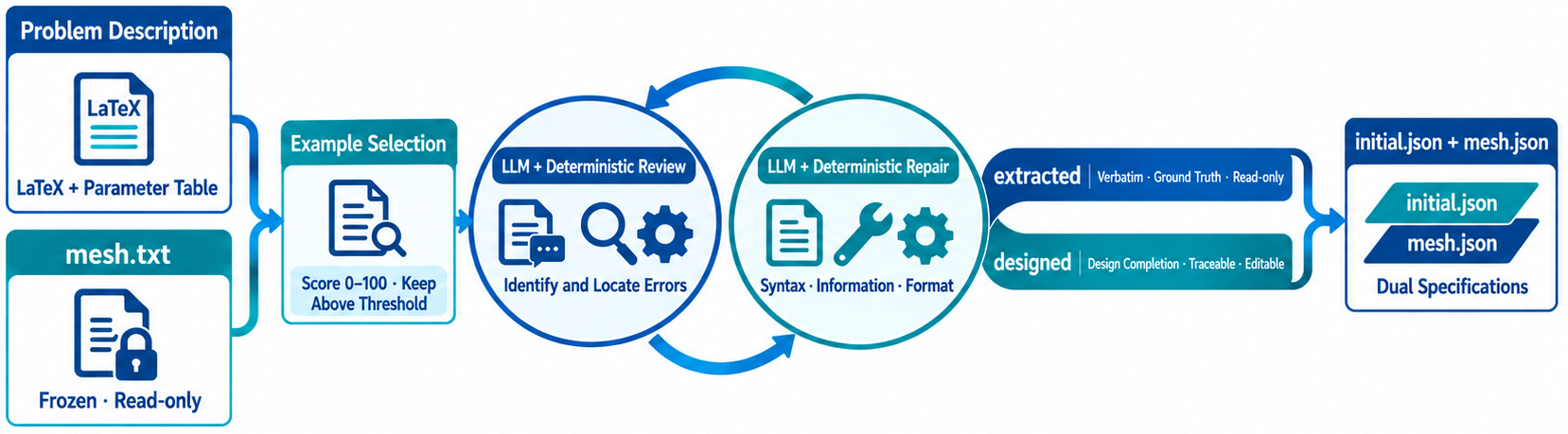}
\caption{The Initialization stage: the description becomes two
provenance-annotated specifications through generation, checking,
review, and targeted repair.}
\label{fig:initialization}
\end{figure}

The generated specifications are accepted only after verification, as shown in Figure~\ref{fig:initialization}. A deterministic
layer checks their syntactic well-formedness, the internal consistency
of labels, roles, and lists, the mutual consistency between the two
specifications, and compliance with a number of physical plausibility rules, and normalizes
formatting irregularities mechanically. A second layer, the semantic
reviewer, compares the specifications against the source document: nothing stated in the
source may be omitted or contradicted, and nothing designed may be
physically implausible. Both layers run every round. If either
layer reports an issue, a repair step rewrites only the flagged entries and
leaves everything else untouched, after which the loop re-checks. The loop
repeats until both layers pass or a given repair budget is exhausted;
when the budget is exhausted, the case is set aside for human
review rather than being allowed to continue with known defects. A provided mesh script,
when one exists, is treated as immutable input: it is checked once,
frozen, and all mesh metadata is derived from it.
Table~\ref{tab:initchecks} summarizes the items covered by the two
layers and how each one is verified. 

\begin{table}[!t]
\caption{The items checked on the generated specifications in the
Initialization stage.}
\label{tab:initchecks}
\centering
\small
\begin{tabular}{p{3.8cm}p{8.4cm}p{3.9cm}}
\toprule
\textbf{Item} & \textbf{Contents} & \textbf{How} \\
\midrule
Provenance annotations & extracted/designed provenance, design
rationale & Deterministic \\
Extraction completeness & constants, boundary conditions, objectives, requirements & Semantic review \\
Transcription fidelity & values, formulas, weak-form
signs & Semantic review \\
Boundary and label rules & disjoint lists, roles, known labels & Deterministic \\
Field declarations & types, approximation spaces, test-function match & Deterministic + semantic \\
Units & canonical spelling, dimensional consistency &
Deterministic + semantic \\
Initial design & parseability, solid presence, connectivity, prescribed seed, solid fraction &
Deterministic + semantic \\
Mesh definition & mesh construction, border labels, domain dimensions, comments, identifiers &
Deterministic \\
Cross-file consistency & labels, edges, lists &
Deterministic \\
Physical plausibility & constants, settings & Semantic review \\
Review consistency & consistent verdicts across rounds & Semantic
review \\
\bottomrule
\end{tabular}
\end{table}

\subsection{Generation}

The second stage produces the solver
program. In addition to the results of the Initialization stage, the
writer receives the reference programs and meshes selected during
that stage: the selected main programs exemplify the usage patterns
of the language, while the selected meshes illustrate mesh-construction
conventions.
The writer translates the
specification into the AutoFreeFEM symbolic language: the domain, the
fields and test functions, the material constants, the material
contrast, the objective, and the weak forms are expressed symbolically,
and the level-set method is configured with the mesh, the
mesh-size bounds, the advection settings, the iteration
count, the initial $\phi$, and the boundary-label lists.
The resulting program provides a structured representation from which the
AutoFreeFEM library later generates the FreeFem++ scripts together
with \LaTeX{} documentation of the problem.

\begin{figure}[!tb]
\centering
\includegraphics[width=0.9\textwidth]{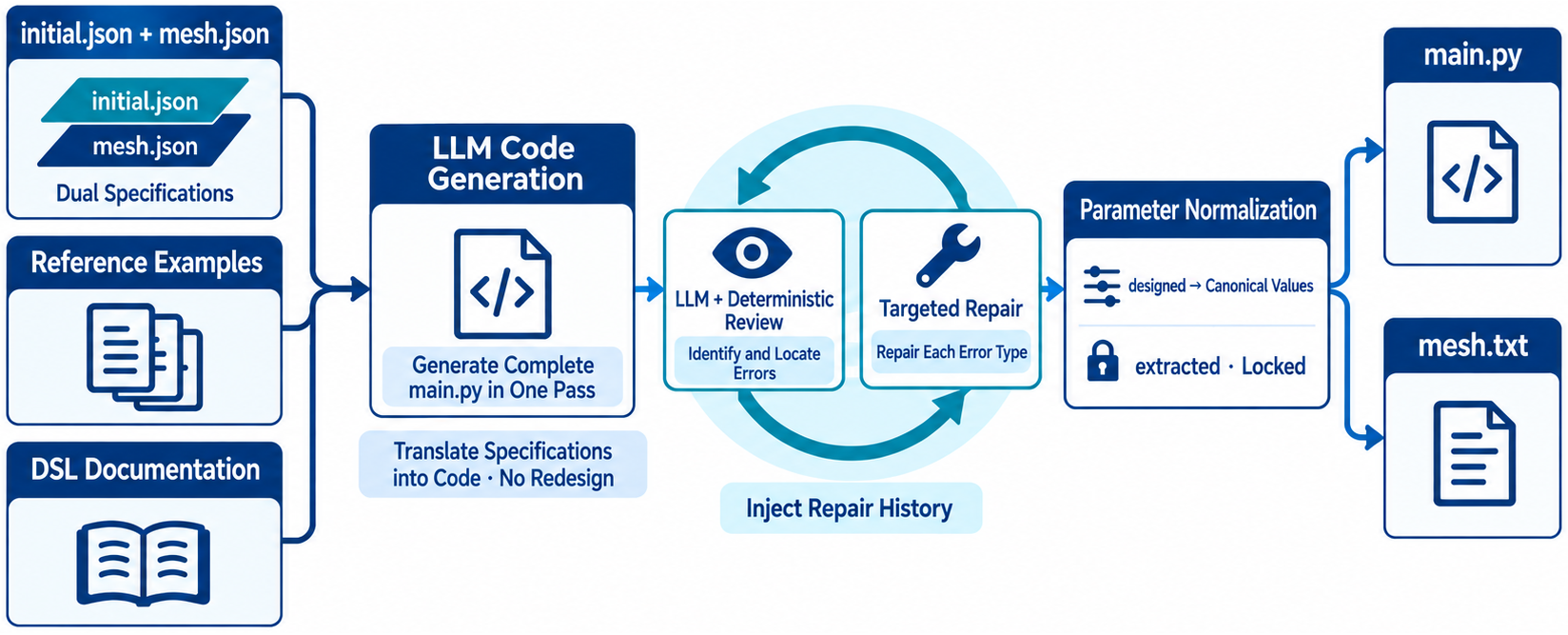}
\caption{The Generation stage: the solver program is written from the
specification with the selected reference examples, and the static
repair and normalization steps finalize it.}
\label{fig:generation}
\end{figure}

The draft then passes through a static repair loop mirroring the first
stage, with the same repair budget. A set of deterministic checks examines the program
structure, the call forms, the transcription of the specification, and
the validity of the parameter settings; a
semantic reviewer checks, against the specification and the
description, that every formula and weak form is mathematically
equivalent to the specified one and that the program uses no feature
the language cannot express. Together, the two layers aim to ensure
that the generated solver program is structurally valid and consistent
with the mathematical specification. Table~\ref{tab:genchecks}
summarizes what is checked and normalized at this stage.
Figure~\ref{fig:generation}
illustrates the stage.

\begin{table}[!t]
\caption{The items checked and normalized on the generated program in
the Generation stage.}
\label{tab:genchecks}
\centering
\small
\begin{tabular}{p{4.5cm}p{8.6cm}p{3.0cm}}
\toprule
\textbf{Item} & \textbf{Contents} & \textbf{How} \\
\midrule
Program structure & import order, domain declaration, problem name, comments & Deterministic \\
Call forms & argument counts, names, types &
Deterministic \\
Transcription & labels, initial design, iteration count, mesh filename, extracted
constants & Deterministic \\
Expressible weak forms & component access, boundary-term
labels, vector gradients, integral division, Dirichlet collisions &
Deterministic \\
Penalty rules & placement, range, notation &
Deterministic \\
Mathematical equivalence & formulas, weak forms & Semantic
review \\
Parameter normalization & material-cost weight, mesh
size, diffusion coefficient and advection step, initial design &
Deterministic \\
\bottomrule
\end{tabular}
\end{table}

\subsection{Execution}

The pipeline shown in Figure~\ref{fig:executionreflection} launches the generated program, which
emits the FreeFem++ scripts, and then runs the optimization. Because
FreeFem++ does not report completion itself, the runner watches the
stream of iteration snapshots the optimizer writes: a new snapshot
resets a per-step clock, and a step that stalls for too long is
terminated, so no run remains stalled indefinitely; the reason
for any termination, whether a solver error, a stalled step, or a completed
design loop, is recorded with the run.

A parameter sweep is then
launched automatically. The material-cost weight, the
diffusion coefficient, and the advection step are varied over a fixed parameter grid. Each
combination runs in its own copy of the program in a bounded worker
pool; within each run, the step attaining the smallest objective value
over that run's history is taken as the run's result. For every parameter combination, this step produces
a complete record of the code, parameters, exit status, objective
history, and final design; the final design reported for each case is then selected by human inspection from these records.

\subsection{Reflection}

\begin{figure}[!tb]
\centering
\includegraphics[width=0.95\textwidth]{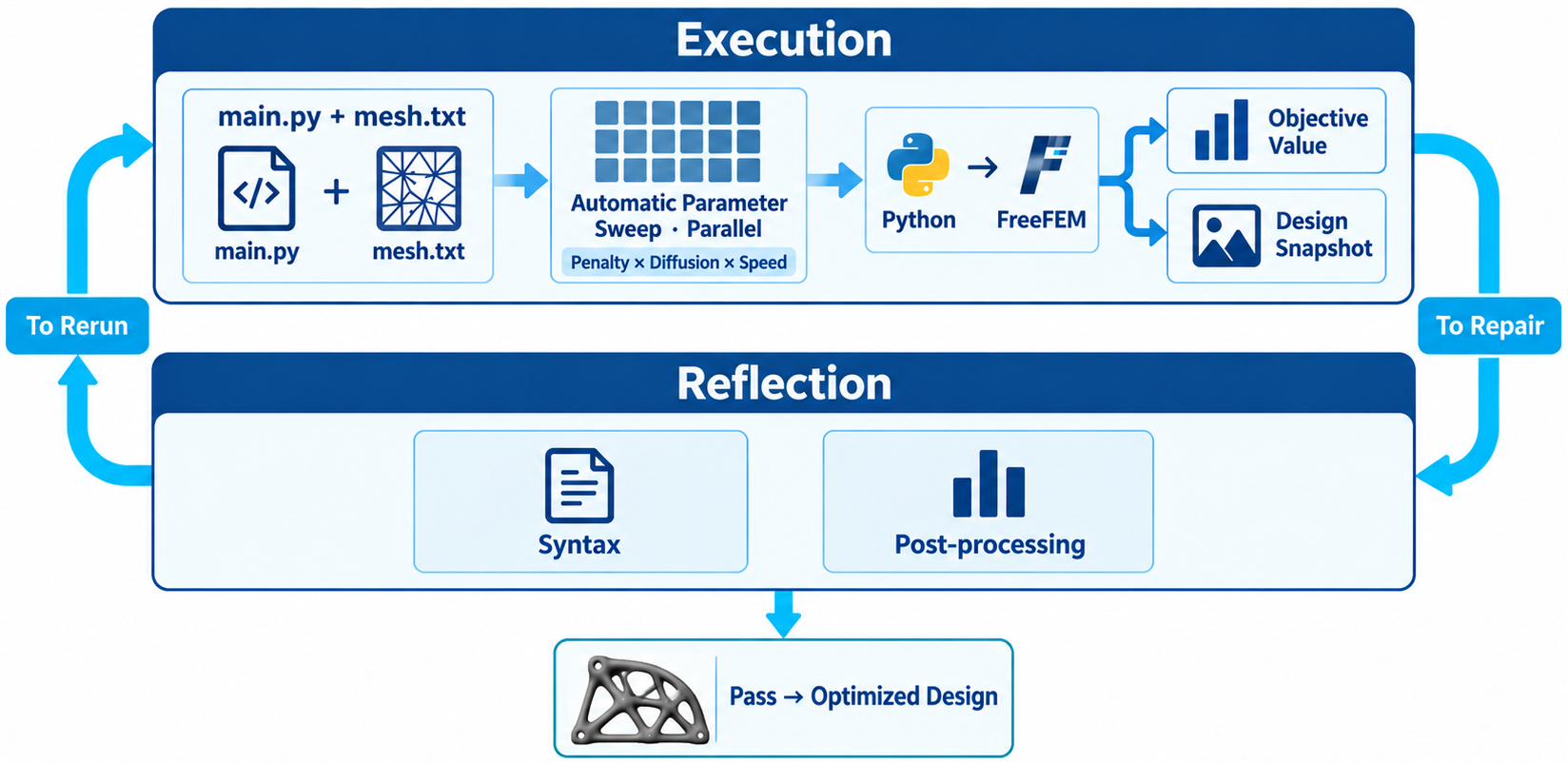}
\caption{The Execution and Reflection stages: parameter-swept runs
that produce no design snapshot are classified as program errors,
repaired, and re-executed; a completed sweep passes the overall design
check and is post-processed into the final figures.}
\label{fig:executionreflection}
\end{figure}

The last stage shown in Figure~\ref{fig:executionreflection} closes the loop between execution
and repair, and serves as the final validation stage: the static checks and the
reviews of the initialization and generation stages already intercept
the vast majority of errors, so reflection performs no repair in
most cases. Reflection has two duties. First, when the executed
group produces no design snapshot at all, the failure is classified as
a program error: the program failed before writing any design. The
text model then examines the program and repairs the code, and the
repaired code is written back into the same group and re-executed; a
program repair never creates a new group. Second, once the executed group
produces design snapshots, reflection first performs an overall check:
the vision model inspects the final design of every experiment against
the problem and assesses whether at least one design is physically
plausible. When
the check fails, a design repairer, responsible for formulas,
parameters, and other design-related issues, proposes a new group of
variants, each a minimal edit of the code that already ran; a
deterministic sanitizer enforces the execution contract on every
variant; and the new group is executed and checked again, closing the
execution--reflection loop. The loop repeats until the overall check
passes or the repair budget is exhausted, and each round is packaged
into the growing history together with its run evidence and images.
The checked designs are then post-processed: their screenshots are
composed into the figure mosaics and the temporary files are cleaned
up. Table~\ref{tab:refchecks} summarizes the
examinations applied in this stage.

\begin{table}[!t]
\caption{The examinations applied in the Execution and Reflection
stages.}
\label{tab:refchecks}
\centering
\small
\begin{tabular}{p{3.7cm}p{8.8cm}p{3.6cm}}
\toprule
\textbf{Item} & \textbf{Contents} & \textbf{How} \\
\midrule
Execution contract & solver invocation, solver search, script writing & Deterministic \\
Program repair trigger & no design snapshot in the executed group
$\Rightarrow$ program error & Deterministic \\
Program scope & program writing, mesh construction & Text-model
detector \\
Overall check & at least one correct design among the experiments & Vision-model
judgment \\
Post-processing & screenshot mosaics, temporary-file cleanup & Deterministic \\
History packaging & code, parameters, consoles, images & Deterministic \\
\bottomrule
\end{tabular}
\end{table}

% =====================================================================
\section{Results}

This section reports the experimental evaluation of Agent4FreeFem.
All experiments used the DeepSeek-V4-Pro and DeepSeek-V4.1-Flash
models with deterministic decoding. The Flash model performed the
overall check of the Reflection stage, which consumes design images,
while all other stages ran on the Pro model. The experiments were
conducted on a single workstation equipped with an AMD Ryzen~7
9800X3D CPU, an
NVIDIA GeForce RTX~5070 Ti GPU, and 64\,GB of memory. The workstation
ran Windows~11 with WSL2, under which the FreeFem++ solvers were
executed. We report, in order, the overall statistics across the
suite, eight representative designs covering the linear/nonlinear and
thermal-diffusion/elasticity combinations, a component ablation study,
error-injection experiments, and two extension cases: fluid flows and
thermo-elasticity.

\subsection{Overall Statistics}

Each case follows the same workflow: the initialization
and the generation stages produce one program template from the
problem description, and this template is then executed under each of 64
parameter combinations, using a separate working copy for each run;
Table~\ref{tab:sweepgrid} lists the parameter values swept.
Table~\ref{tab:runtime} reports the number of nodes of the
initial mesh of each case. The meshes consist of unstructured triangular
elements. During optimization,
the meshes are reconstructed to conform to the evolving material
interfaces, resulting in changes in the numbers of mesh vertices and
elements. The execution and reflection stages form a closed loop: a run
that produces no design snapshot at all is classified as a program
error, repaired, and re-executed, while a run that completes first
undergoes the overall design check and is then post-processed.
Whenever a run failed, it did so during
startup, before any substantial computation started, and each such
failure cost one program-repair round followed by a re-execution; the
final round of each case ran the full 64-combination sweep to
completion, passed the overall check, and was post-processed. In all
the cases of the reported experiments, the overall check passed on the
first attempt, no design-repair round was ever needed, which
indirectly attests to the effectiveness of the large language model
and of the overall architecture. The run
counts in
Table~\ref{tab:runtime} report
these loop rounds; the table also reports the runtimes.
The 32 cases involved 2{,}048 completed parameter-sweep
runs. Runtime completion indicates that the numerical optimization
finished successfully, whereas design quality is assessed separately
through the reflection-stage review and subsequent human inspection.
The 64 runs of
a sweep are computed in parallel on a bounded worker pool, so the
execution times reported here are wall-clock times of the parallel
sweep, not the sum of the individual runs.

\begin{table}[!t]
\caption{The parameter values swept by the pipeline. Each case is
executed once under every combination of the three parameters, giving
$4^3=64$ combinations per case; the diffusion coefficient and the
advection step share the same grid for the elasticity and the thermal
diffusion groups.}
\label{tab:sweepgrid}
\centering
\small
\begin{tabular}{l r r r}
\toprule
Group & material-cost weight $\ell$ & diffusion coefficient $\alpha$ & advection step $\tau$ \\
\midrule
Elasticity & $\{1, 2.5, 5, 10\}$ & $\{10^{-2}, 10^{-3}, 10^{-4}, 10^{-5}\}$ & $\{2\times10^{-3}, 2\times10^{-4}, 2\times10^{-5}, 2\times10^{-6}\}$ \\
Thermal diffusion & $\{25, 50, 100, 500\}$ & $\{10^{-2}, 10^{-3}, 10^{-4}, 10^{-5}\}$ & $\{2\times10^{-3}, 2\times10^{-4}, 2\times10^{-5}, 2\times10^{-6}\}$ \\
\bottomrule
\end{tabular}
\end{table}

The failures encountered across the suite were handled by the
automated checking and repair loops rather than by a human debugger.
Every repair loop is subject to a budget of one initial generation
followed by at most ten repairs,
eleven rounds in total, and a case that fails to pass within the
cap is set aside for human review instead of proceeding. The
maximum number of rounds in Table~\ref{tab:runtime} was 10 (less than 11), so no
case ever exhausted its budget. In
addition, the final results of all cases were inspected by human
experts, providing additional evidence of the numerical validity of
the solutions.

These numbers substantiate the properties claimed in the
Introduction. \textbf{Efficient connectivity}: every case ran from the
description to the finished sweep without any manual handoff between
stages, at a mean end-to-end runtime of 85 minutes 2 seconds per case,
so the reported time per case is the complete cost of the automated
workflow, from text to the optimized designs. \textbf{Process autonomy}: the specification, the
reference-example selection, the code generation, and the parameter
sweep of every case were completed autonomously by the pipeline, with
no solver code written or debugged by hand. \textbf{Process
verifiability}: every failure was caught by the deterministic checks,
the semantic review, or the execution evidence, then classified,
repaired, and re-executed until it passed; the component ablation
study of Section~4.3 further quantified the contribution of the
components to the success rate and the repair overhead.

\begin{table}[!t]
\caption{Per-case runtime statistics over the 32-case suite, compiled
from the pipeline's per-case runtime summary. The Nodes column lists
the number of mesh nodes of the initial mesh. The initialization and
generation columns report per-case totals; the execution and reflection column reports
the wall-clock time of the 64 parallel runs plus the reflection
time. Each initialization and generation run count comprises one
initial generation plus at most ten repairs; the execution and
reflection run count is the number of loop rounds of Section~4.1. }
\label{tab:runtime}
\centering
\footnotesize
\resizebox{\textwidth}{!}{%
\begin{tabular}{l l r r r r r r r r}
\toprule
Type & Case & \multicolumn{3}{c}{Initialization} &
\multicolumn{2}{c}{Generation} &
\multicolumn{2}{c}{Execution + Reflection} & Total \\
\cmidrule(lr){3-5} \cmidrule(lr){6-7} \cmidrule(lr){8-9}
& & Nodes & Runs & Time & Runs & Time & Runs & Time & \\
\midrule
\multirow{8}{*}{Linear elasticity}
& bridgePier & 6{,}932 & 3 & 17m08s & 3 & 26m49s & 3 & 76m02s & 119m59s \\
& cantileverBeam & 2{,}947 & 5 & 30m31s & 3 & 13m42s & 1 & 18m56s & 63m09s \\
& elasticFoundation & 3{,}169 & 2 & 10m20s & 5 & 38m47s & 3 & 71m21s & 120m28s \\
& gravityDam & 2{,}947 & 6 & 35m45s & 1 & 9m26s & 2 & 41m11s & 86m22s \\
& hangingBracket & 2{,}947 & 2 & 11m26s & 3 & 38m06s & 1 & 10m34s & 60m06s \\
& mbbBeam & 4{,}489 & 3 & 22m17s & 3 & 44m12s & 1 & 17m37s & 84m06s \\
& pressureVessel & 2{,}947 & 3 & 20m32s & 4 & 34m34s & 2 & 40m26s & 95m32s \\
& shearWall & 2{,}947 & 5 & 33m45s & 3 & 17m09s & 1 & 8m11s & 59m05s \\
\midrule
\multirow{8}{*}{Nonlinear elasticity}
& bridgePierNonlinear & 6{,}932 & 5 & 32m11s & 1 & 5m05s & 1 & 39m36s & 76m52s \\
& cantileverBeamNonlinear & 2{,}947 & 10 & 59m14s & 3 & 25m21s & 1 & 77m08s & 161m43s \\
& elasticFoundationNonlinear & 3{,}169 & 4 & 25m20s & 8 & 68m26s & 1 & 34m10s & 127m56s \\
& gravityDamNonlinear & 2{,}947 & 3 & 37m26s & 5 & 57m51s & 1 & 53m05s & 148m22s \\
& hangingBracketNonlinear & 2{,}947 & 4 & 30m50s & 3 & 27m59s & 1 & 50m16s & 109m05s \\
& mbbBeamNonlinear & 4{,}489 & 9 & 114m52s & 4 & 35m14s & 2 & 140m05s & 290m11s \\
& pressureVesselNonlinear & 2{,}947 & 2 & 11m00s & 3 & 21m15s & 1 & 44m36s & 76m51s \\
& shearWallNonlinear & 2{,}947 & 2 & 12m08s & 1 & 3m58s & 1 & 35m25s & 51m31s \\
\midrule
\multirow{8}{*}{Linear thermal diffusion}
& contactPadCooling & 4{,}373 & 5 & 29m46s & 1 & 2m51s & 2 & 26m40s & 59m17s \\
& doubleSidedCooling & 3{,}721 & 8 & 36m20s & 3 & 21m32s & 3 & 42m10s & 100m02s \\
& doubleSinkColdPlate & 22{,}656 & 2 & 46m32s & 1 & 8m18s & 2 & 24m05s & 78m55s \\
& finCrossSection & 3{,}721 & 2 & 7m15s & 1 & 3m46s & 1 & 6m17s & 17m18s \\
& heatedPanel & 3{,}721 & 6 & 39m31s & 1 & 2m44s & 1 & 2m46s & 45m01s \\
& microchannelCell & 3{,}721 & 5 & 47m55s & 3 & 14m18s & 2 & 26m35s & 88m48s \\
& timPartialContact & 4{,}072 & 5 & 29m16s & 1 & 6m25s & 1 & 3m34s & 39m15s \\
& vtpColdPlate & 10{,}134 & 4 & 30m54s & 1 & 3m50s & 1 & 7m37s & 42m21s \\
\midrule
\multirow{8}{*}{Nonlinear thermal diffusion}
& contactPadCoolingNonlinear & 4{,}010 & 2 & 13m03s & 4 & 26m03s & 1 & 11m45s & 50m51s \\
& doubleSidedCoolingNonlinear & 3{,}721 & 3 & 27m27s & 1 & 4m27s & 1 & 11m52s & 43m46s \\
& doubleSinkColdPlateNonlinear & 3{,}833 & 6 & 37m27s & 1 & 4m43s & 1 & 7m18s & 49m28s \\
& finCrossSectionNonlinear & 3{,}721 & 3 & 17m50s & 2 & 13m45s & 1 & 10m37s & 42m12s \\
& heatedPanelNonlinear & 3{,}721 & 4 & 24m47s & 2 & 10m29s & 2 & 24m49s & 60m05s \\
& microchannelCellNonlinear & 4{,}405 & 3 & 18m27s & 2 & 10m11s & 2 & 35m12s & 63m50s \\
& timPartialContactNonlinear & 4{,}072 & 7 & 44m04s & 2 & 15m06s & 2 & 33m19s & 92m29s \\
& vtpColdPlateNonlinear & 6{,}694 & 7 & 68m18s & 5 & 34m19s & 1 & 13m21s & 115m58s \\
\midrule
& \textbf{Total} & & \textbf{140} & \textbf{1023m37s} & \textbf{84} &
\textbf{650m41s} & \textbf{47} & \textbf{1046m36s} & \textbf{2720m54s} \\
& \textbf{Avg} & & \textbf{4.38} & \textbf{31m59s} & \textbf{2.63} &
\textbf{20m20s} & \textbf{1.47} & \textbf{32m42s} & \textbf{85m02s} \\
\bottomrule
\end{tabular}}
\end{table}

\subsection{Representative Designs: Elasticity and Thermal Diffusion}

The eight cases below are grouped into four elasticity cases and
four thermal diffusion cases, covering the linear variant of
each physics and its corresponding nonlinear variant. In all cases, the void material is a thousand times weaker, i.e., $\chi_0=10^{-3}$. For each
variant, the parameter combination judged most reasonable after human inspection of the sweep results is reported in
Table~\ref{tab:params}.

Each case started from its natural-language
description: the pipeline produced the specification, wrote the solver
program, and swept the 64 parameter combinations. Because the sweep
enumerated the full Cartesian product of the three parameters, many of
the $64=4^3$ combinations yielded degenerate designs that barely
departed from the initial checkerboard, developed disconnected or
fragmented material layouts, or collapsed toward an almost empty or
almost full design; the purpose of the
sweep was to run them all and identify the combinations that yielded
satisfactory designs, and the reported design of each case was the
one judged most reasonable after human inspection of the sweep results.

\begin{figure}[!tb]
\centering
\subfloat{%
  \includegraphics[width=0.42\columnwidth]{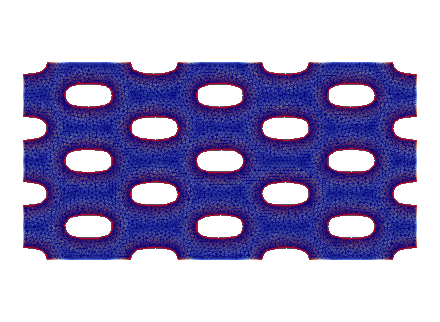}}\hspace{1.5cm}
\subfloat{%
  \includegraphics[width=0.37\columnwidth]{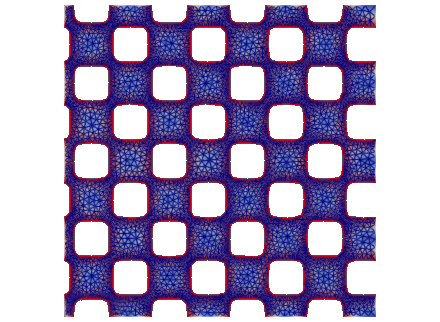}}
\caption{The checkerboard initial design shared by all
elasticity cases (left) and all thermal diffusion cases (right),
shown on representative meshes.}
\label{fig:initial}
\end{figure}

\begin{figure}[!tb]
\centering
\subfloat[Shear wall initial.]{%
  \begin{minipage}[b]{0.32\textwidth}
  \centering
  \includegraphics[width=\linewidth]{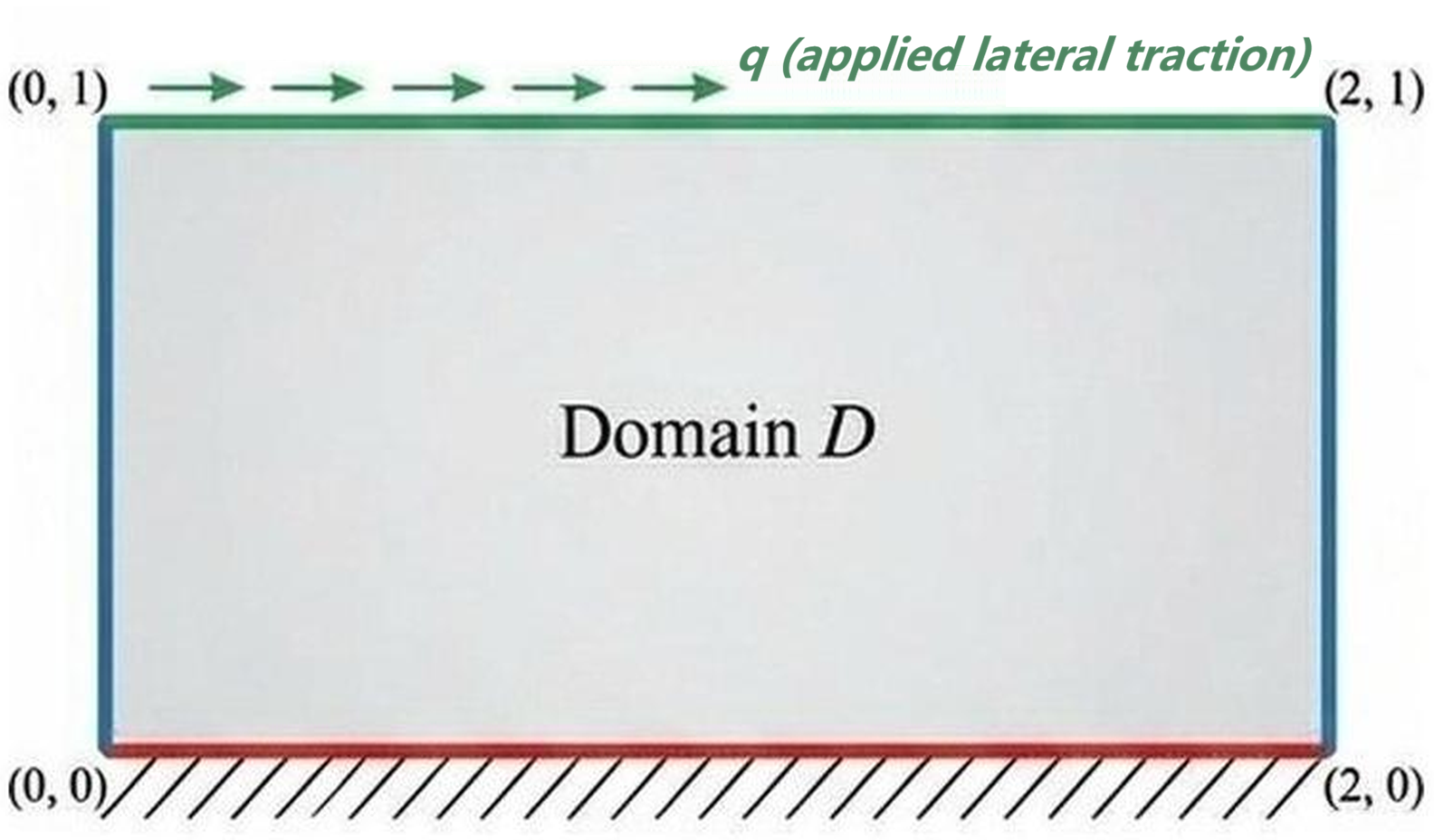}
  \end{minipage}}\hfill
\subfloat[Shear wall, linear result.]{%
  \begin{minipage}[b]{0.32\textwidth}
  \centering
  \includegraphics[width=\linewidth]{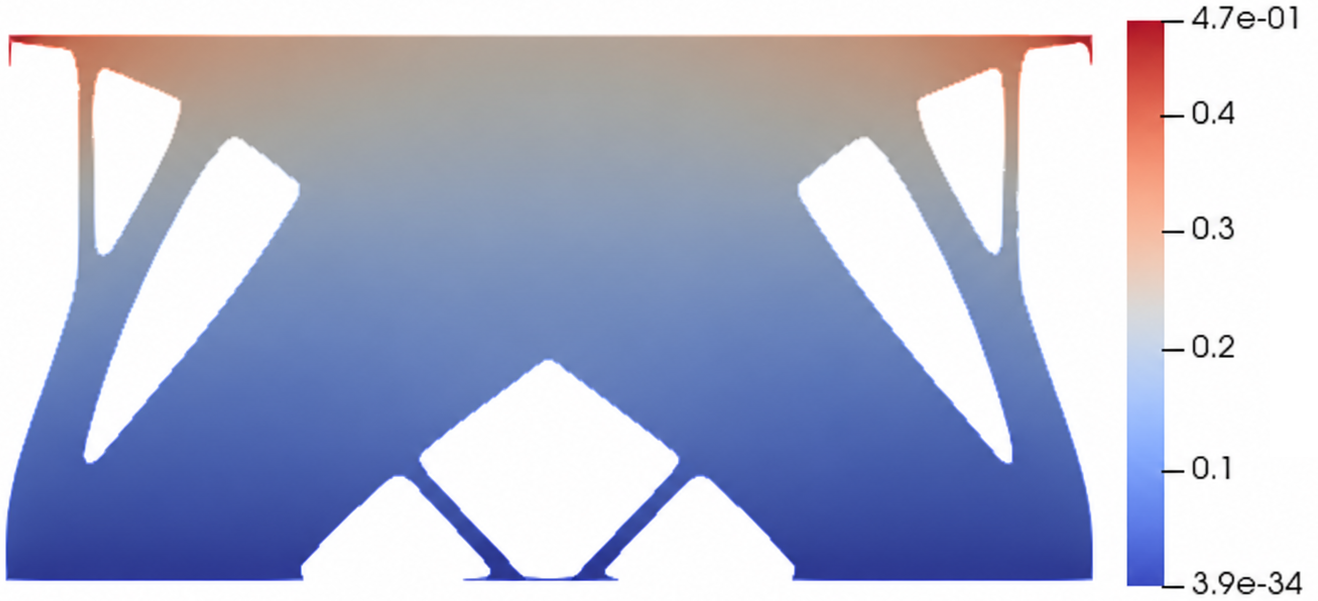}
  \end{minipage}}\hfill
\subfloat[Shear wall, nonlinear result.]{%
  \begin{minipage}[b]{0.32\textwidth}
  \centering
  \includegraphics[width=\linewidth]{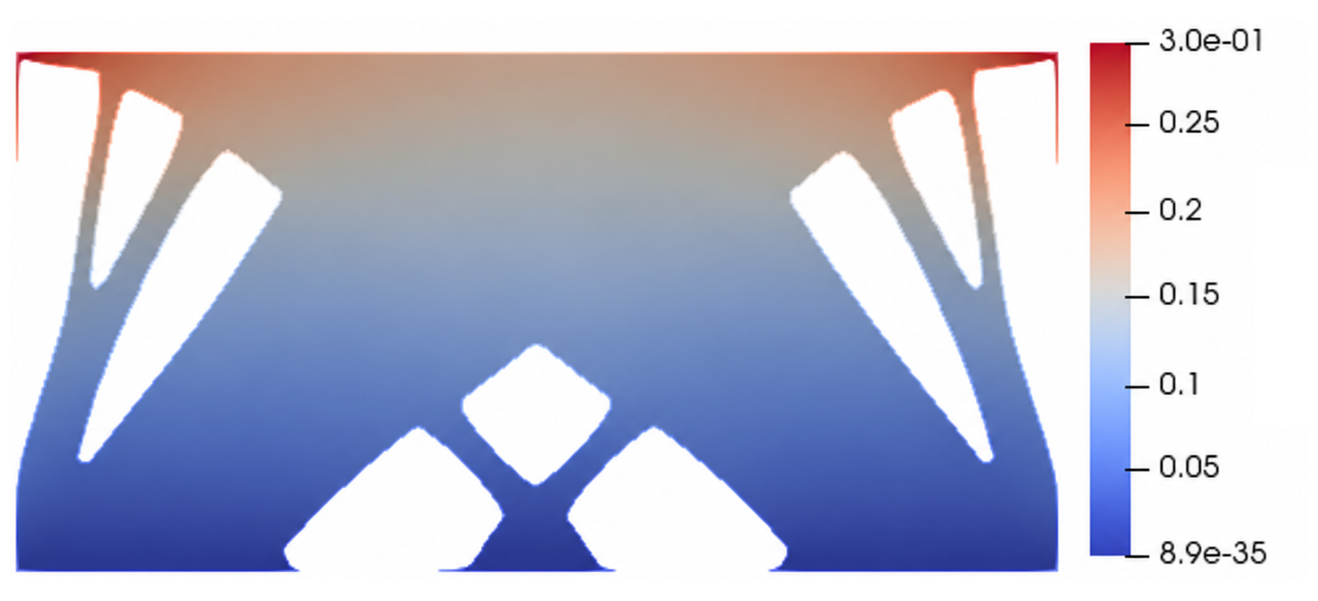}
  \end{minipage}}\\
\subfloat[Pressure vessel initial.]{%
  \begin{minipage}[b]{0.32\textwidth}
  \centering
  \includegraphics[width=\linewidth]{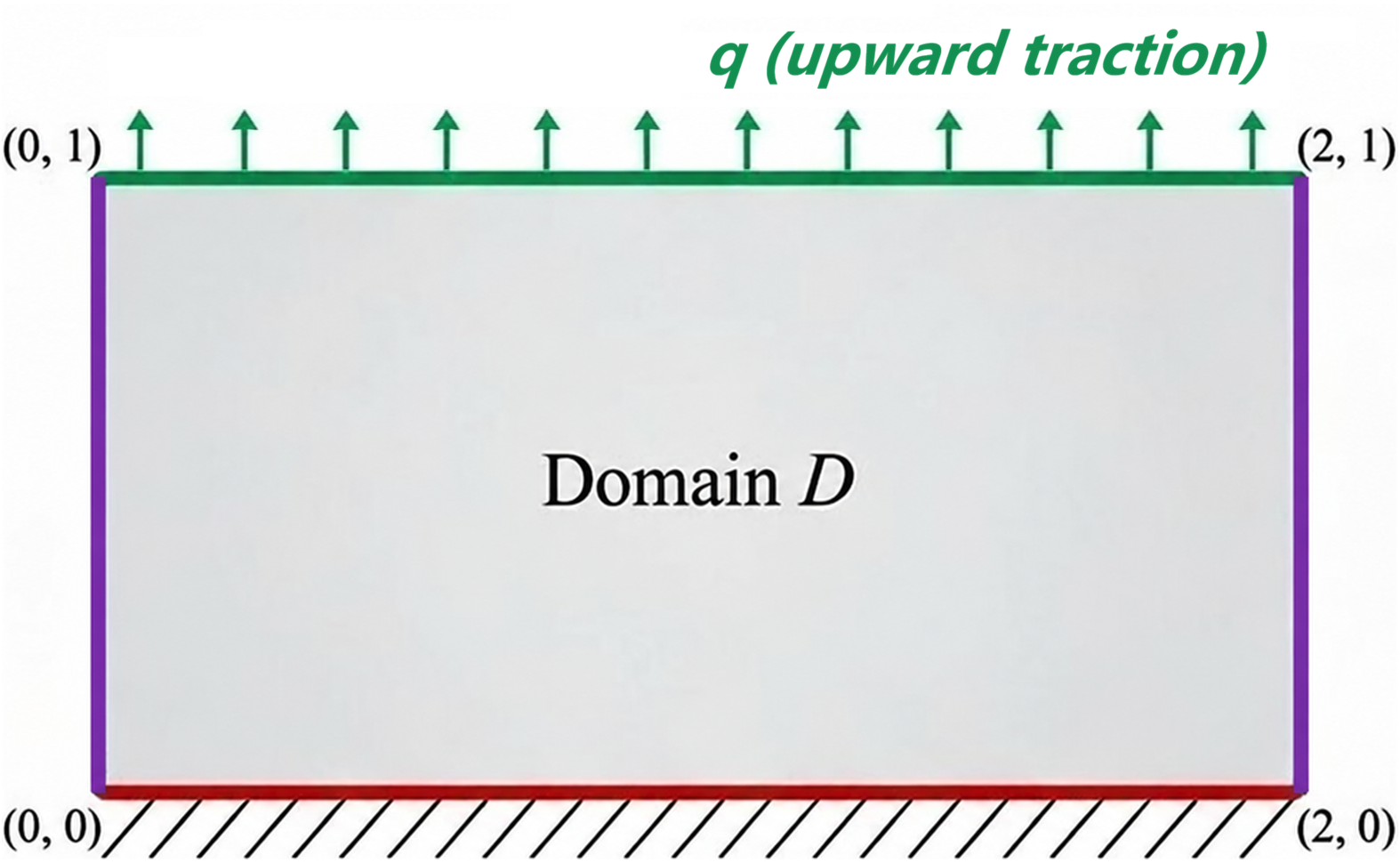}
  \end{minipage}}\hfill
\subfloat[Pressure vessel, linear result.]{%
  \begin{minipage}[b]{0.32\textwidth}
  \centering
  \includegraphics[width=\linewidth]{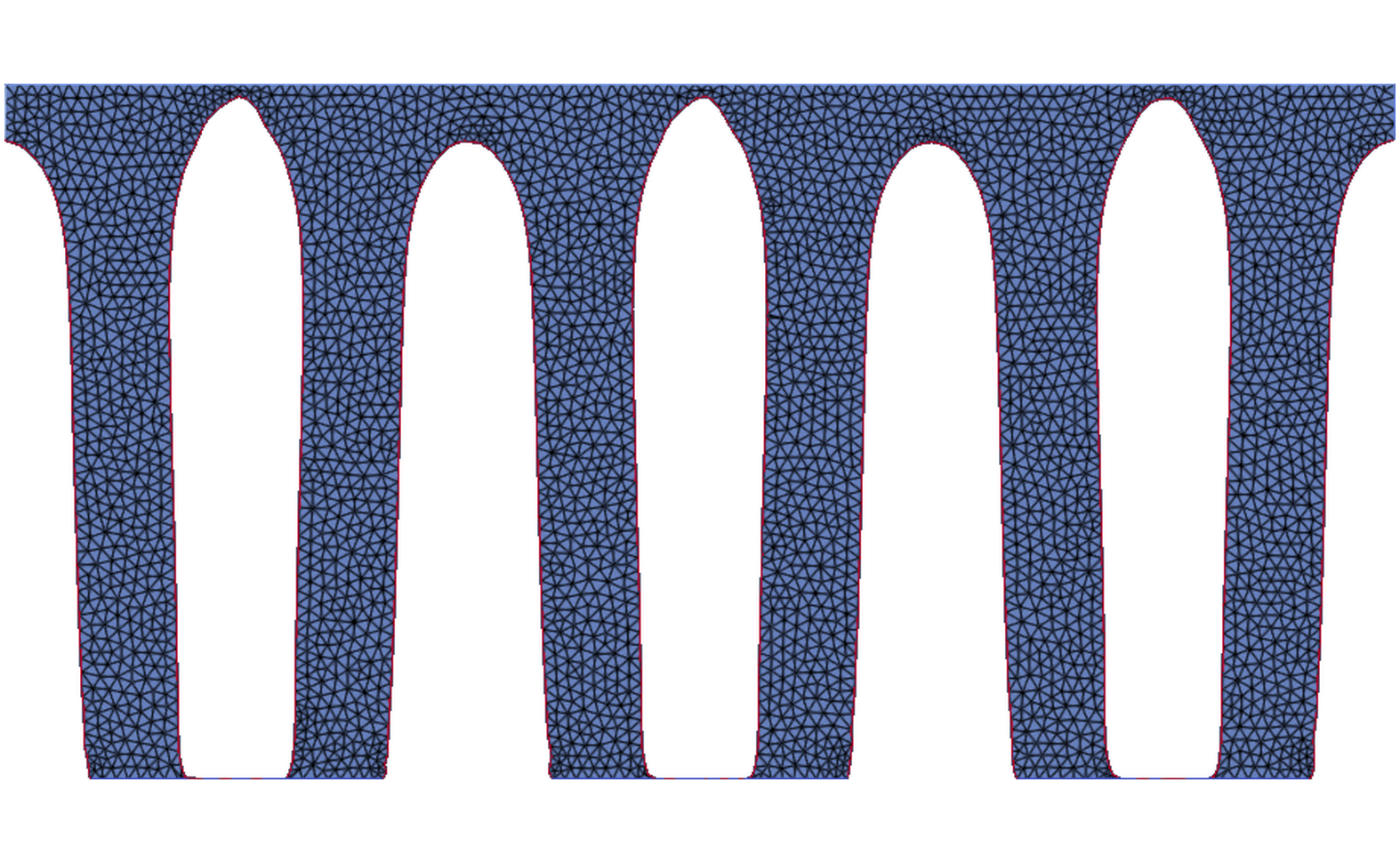}
  \end{minipage}}\hfill
\subfloat[Pressure vessel, nonlinear result.]{%
  \begin{minipage}[b]{0.32\textwidth}
  \centering
  \includegraphics[width=\linewidth]{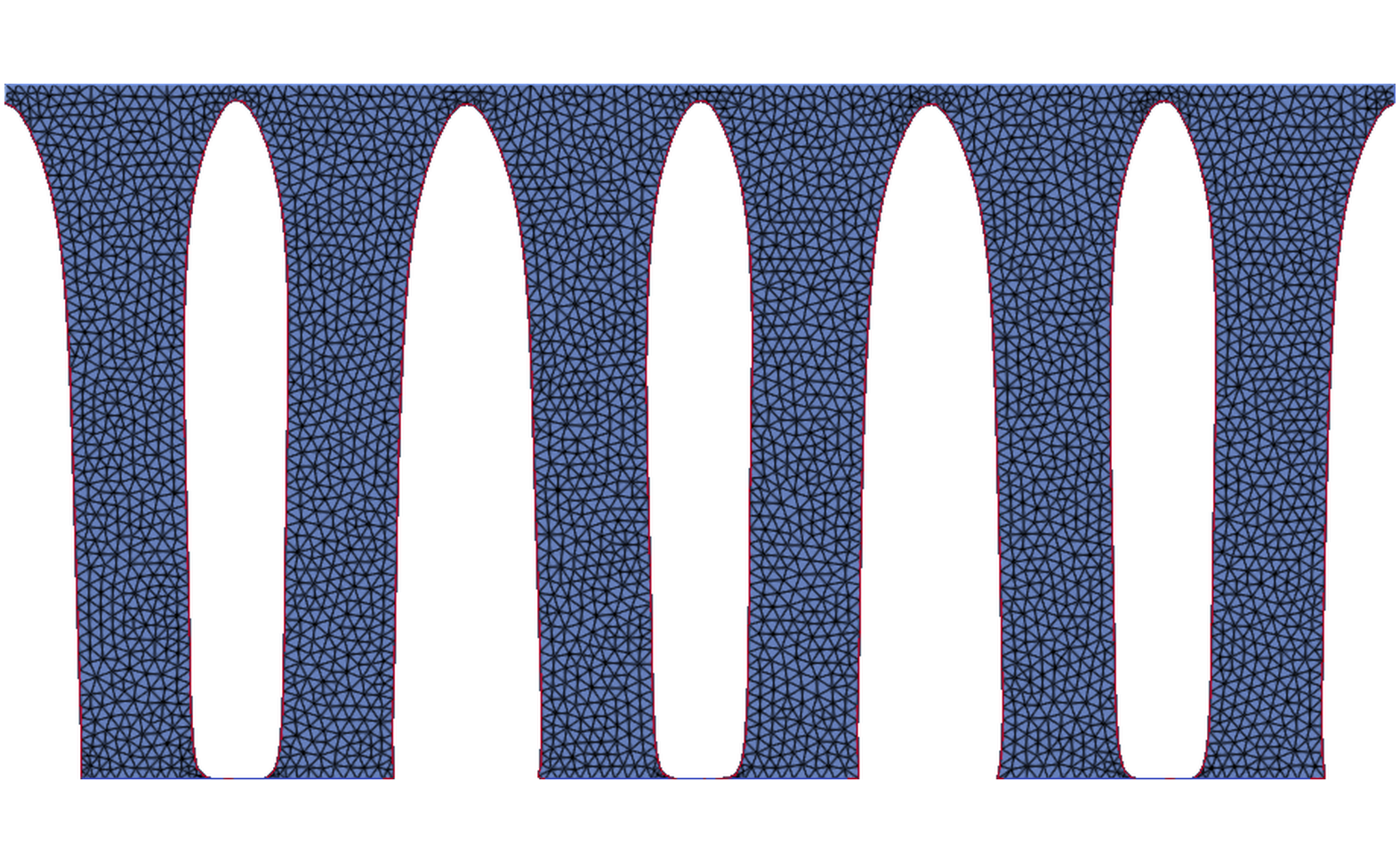}
  \end{minipage}}
\caption{The four elasticity cases: for each problem, the initial
problem schematic (left column) and the final designs selected after human
inspection of the sweep results for the linear (middle column) and nonlinear
(right column) variants; the shear wall panels show the displacement
field, and the pressure vessel panels show the material layout.}
\label{fig:representative}
\end{figure}

\begin{figure}[!tb]
\centering
\subfloat[Cold plate.]{%
  \begin{minipage}[b]{0.32\textwidth}
  \centering
  \includegraphics[width=\linewidth]{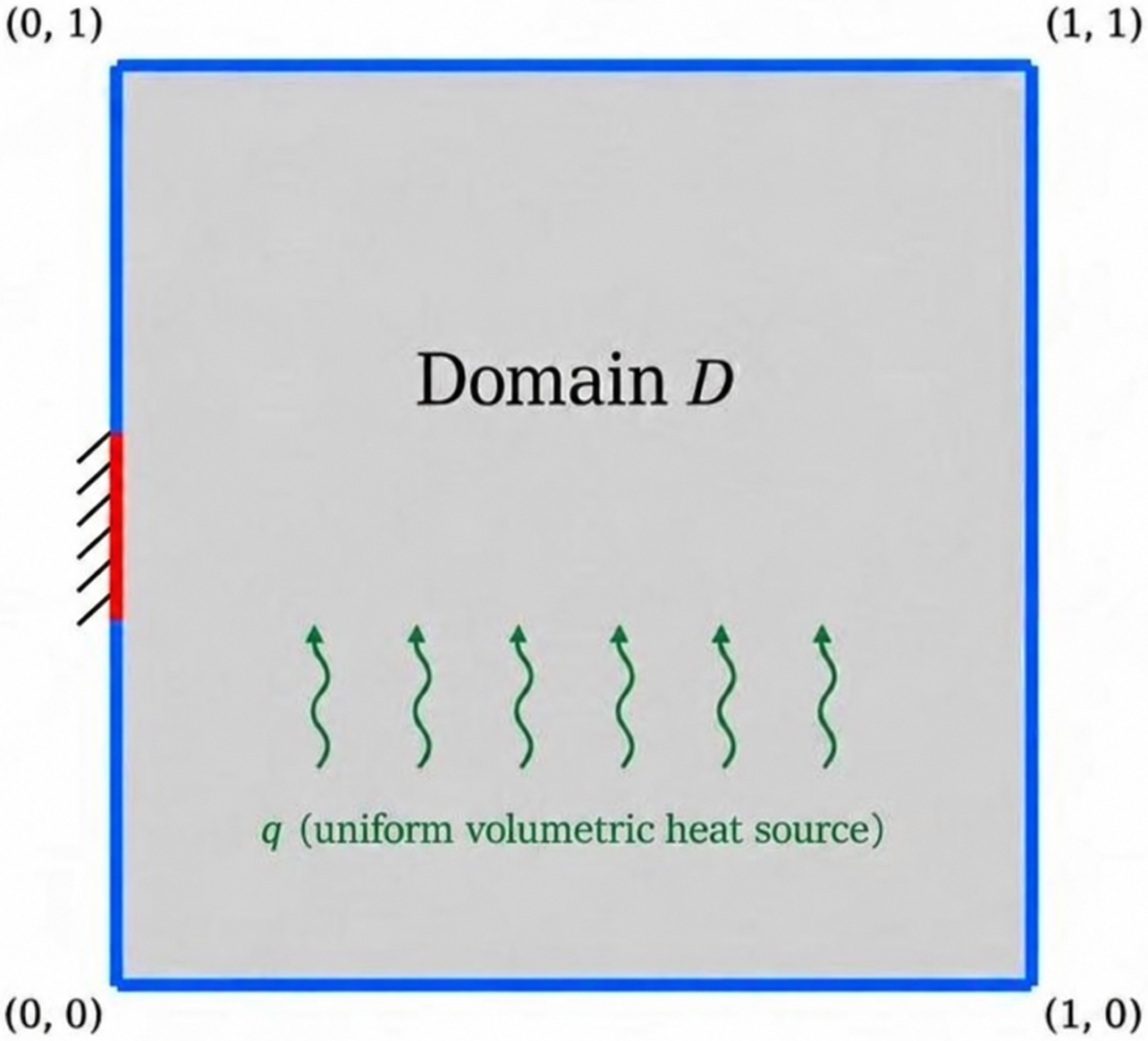}
  \end{minipage}}\hfill
\subfloat[Result on cold plate (linear).]{%
  \begin{minipage}[b]{0.32\textwidth}
  \centering
  \includegraphics[width=\linewidth]{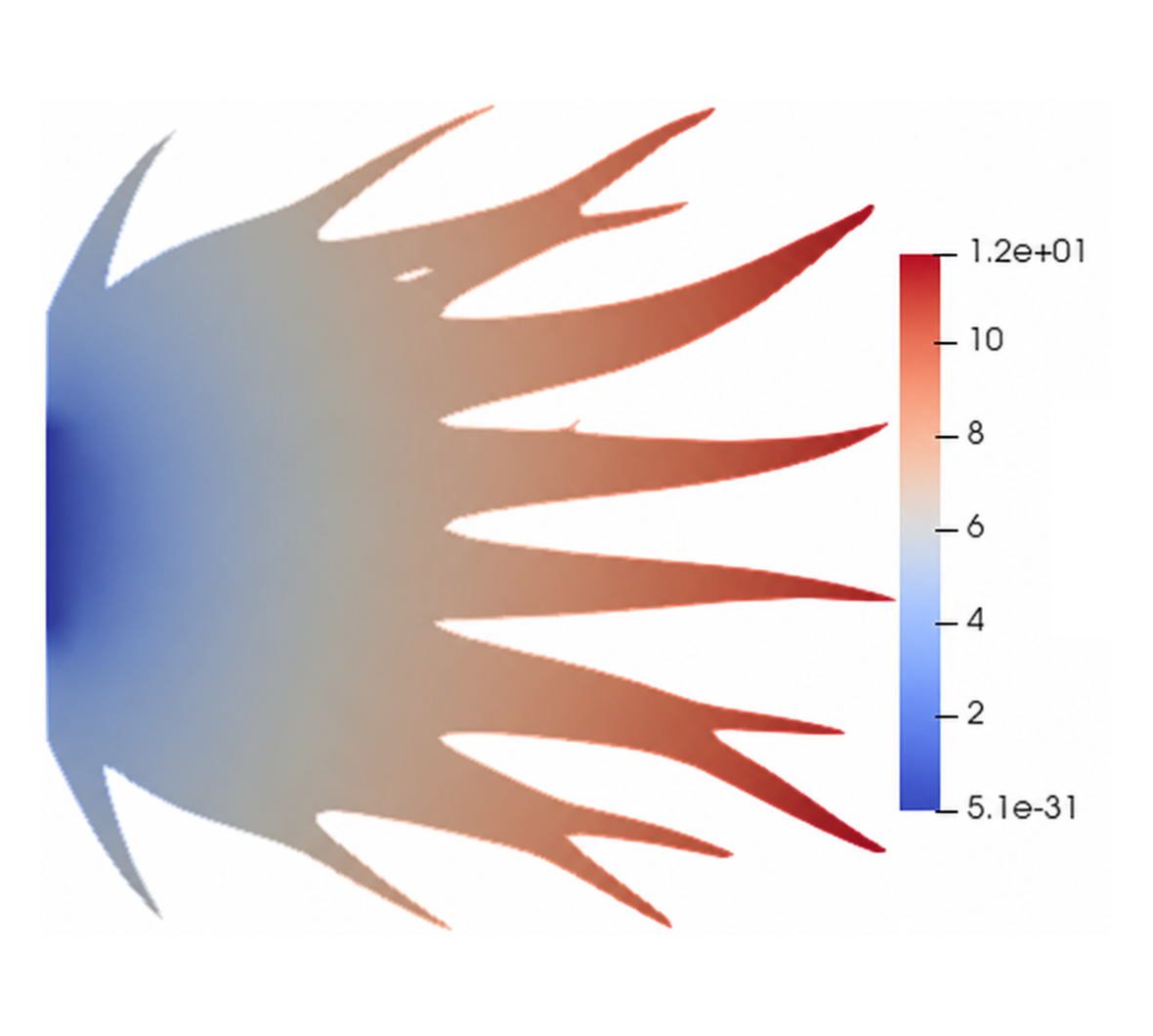}
  \end{minipage}}\hfill
\subfloat[Result on cold plate (nonlinear).]{%
  \begin{minipage}[b]{0.32\textwidth}
  \centering
  \includegraphics[width=\linewidth]{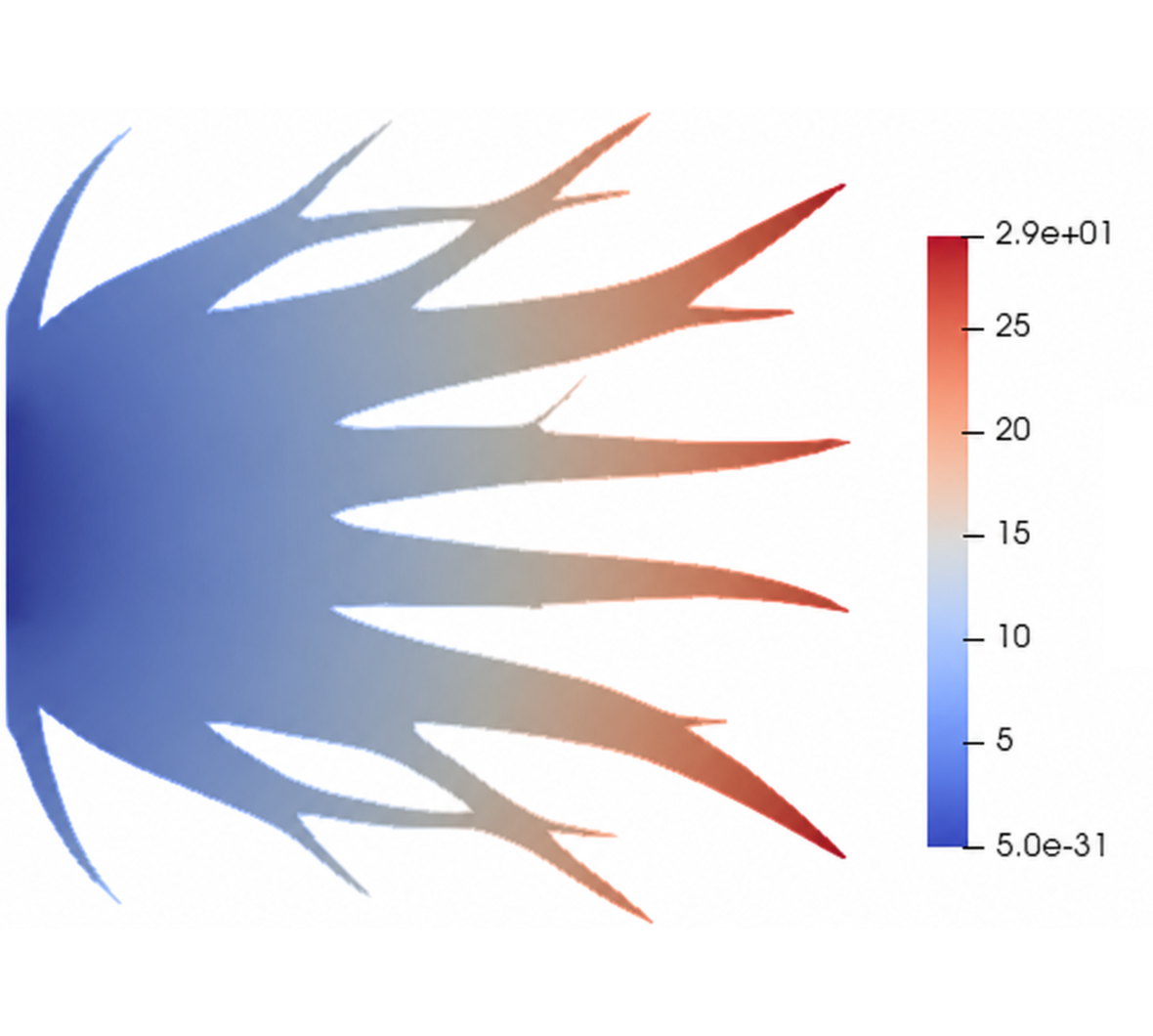}
  \end{minipage}}\\
\subfloat[Two-sink plate.]{%
  \begin{minipage}[b]{0.32\textwidth}
  \centering
  \includegraphics[width=\linewidth]{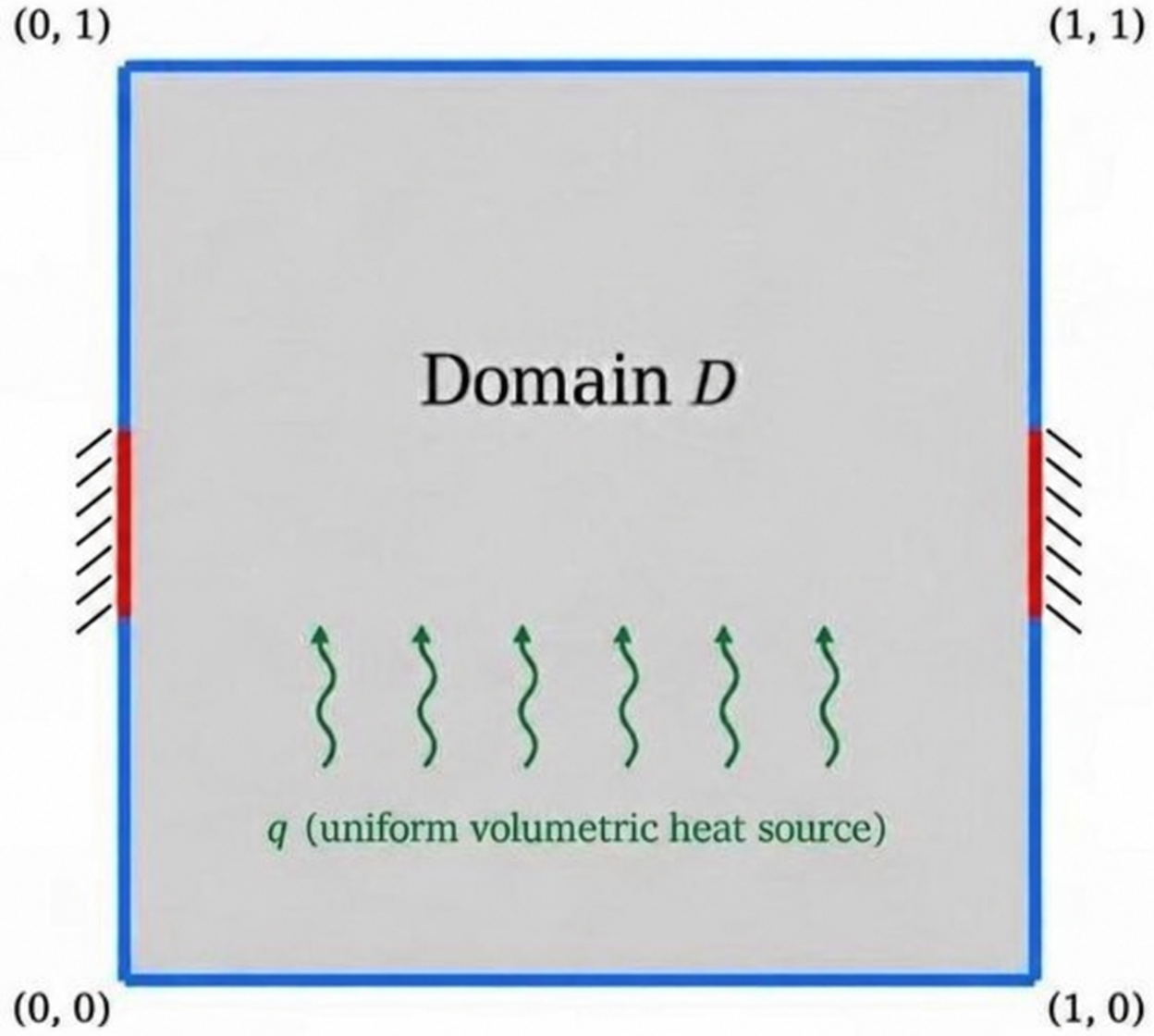}
  \end{minipage}}\hfill
\subfloat[Result on two-sink plate (linear).]{%
  \begin{minipage}[b]{0.32\textwidth}
  \centering
  \includegraphics[width=\linewidth]{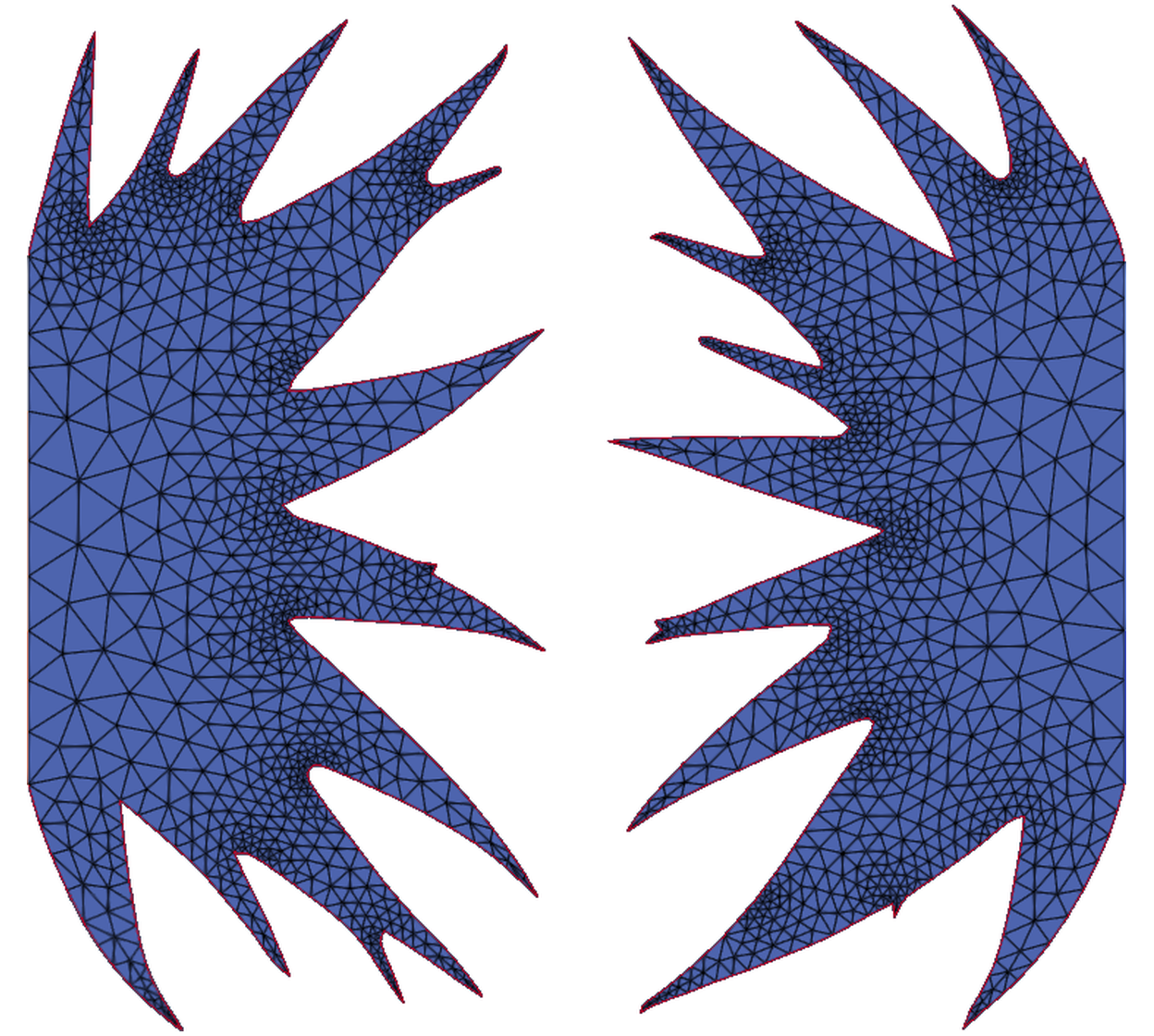}
  \end{minipage}}\hfill
\subfloat[Result on two-sink plate (nonlinear).]{%
  \begin{minipage}[b]{0.32\textwidth}
  \centering
  \includegraphics[width=\linewidth]{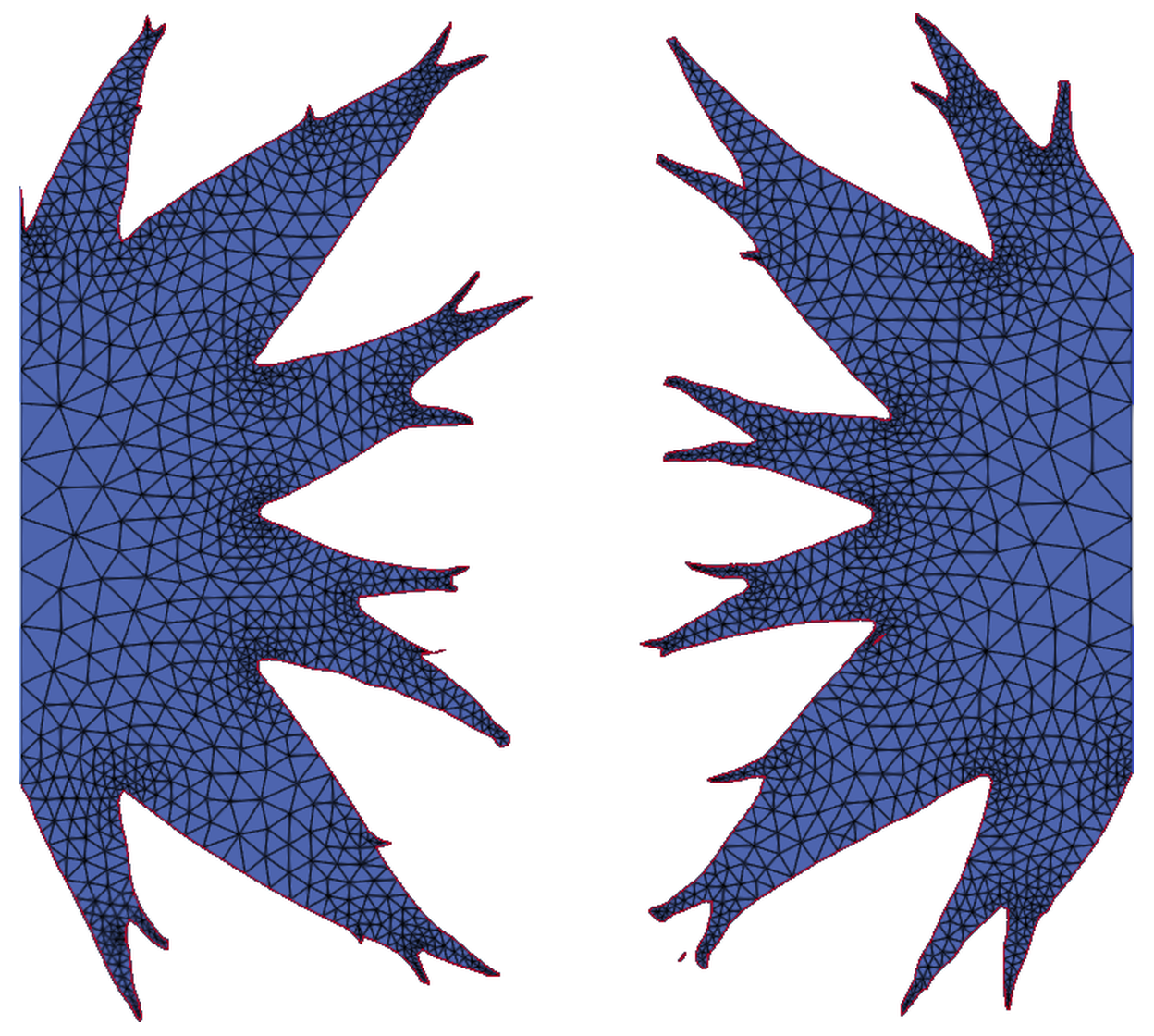}
  \end{minipage}}
\caption{The four thermal diffusion cases: for each problem, the problem schematic (left column) and the final designs selected after human
inspection of the sweep results for the linear (middle column) and nonlinear
(right column) variants; the cold plate panels show the temperature
field.}
\label{fig:diffusion}
\end{figure}

Table~\ref{tab:params} collects these combinations for all eight cases.
Figure~\ref{fig:initial} shows the checkerboard initial designs on the
mesh, and Figures~\ref{fig:representative} and~\ref{fig:diffusion}
show, for the elasticity and the thermal diffusion cases respectively, the
initial problem schematic together with the final designs of the
linear and nonlinear variants.
The two groups also show that problems with and without body
loads are supported alike: the elasticity group carries no body force,
while the thermal diffusion group is driven by body heat.

\begin{table}[!t]
\caption{The parameter combinations judged most reasonable after human
inspection of the sweep results for the eight representative cases.}
\label{tab:params}
\centering
\footnotesize
\begin{tabular}{l c c c c c c c c}
\toprule
& \multicolumn{2}{c}{Shear wall} & \multicolumn{2}{c}{Pressure vessel} &
\multicolumn{2}{c}{Cold plate} & \multicolumn{2}{c}{Two-sink plate} \\
\cmidrule(lr){2-3} \cmidrule(lr){4-5} \cmidrule(lr){6-7} \cmidrule(lr){8-9}
& Linear & Nonlinear & Linear & Nonlinear & Linear & Nonlinear & Linear & Nonlinear \\
\midrule
material-cost weight $\ell$ & 5 & 2.5 & 2.5 & 1 & 100 & 500 & 25 & 25 \\
diffusion coefficient $\alpha$ & $10^{-3}$ & $10^{-3}$ & $10^{-4}$ & $10^{-2}$ & $10^{-4}$ & $10^{-4}$ & $10^{-4}$ & $10^{-4}$ \\
advection step $\tau$ & $2\times10^{-4}$ & $2\times10^{-3}$ & $2\times10^{-5}$ & $2\times10^{-4}$ & $2\times10^{-5}$ & $2\times10^{-6}$ & $2\times10^{-3}$ & $2\times10^{-4}$ \\
\bottomrule
\end{tabular}
\end{table} 

The shear wall is a structural element of a framed building that
receives the entire lateral force of the floors above and carries the
lateral load into the ground. The design domain is the plane-strain
cross-section of the wall, a 2\,m by 1\,m rectangle. The bottom
edge $\Gamma_{\mathrm{fix}}$ is the fixed foundation and the top
edge $\Gamma_{\mathrm{load}}$ carries a uniform tangential
traction, both prescribed; the two lateral edges denoted by $\Gamma_{\rm o}$ are free
for design, as shown in Figure~\ref{fig:representative}(a). The material is linear elastic and the
ersatz constitutive law reads
\begin{equation}\label{fwijeo}
\sigma(\mathbf{u}_s)=\chi(\phi)\,C:\varepsilon(\mathbf{u}_s),
\end{equation}
where $\mathbf{u}_s$ denotes the displacement and $\sigma$ the stress. The objective balances the compliance
term against the material cost. We consider \eqref{lsmTOP} with $\mathcal{J}(\phi)=\int_D
\sigma(\mathbf{u}_s):\varepsilon(\mathbf{u}_s)\,\mathrm{d}x$ subject to
\begin{equation}
\begin{gathered}
%\min_{\phi}\ J(\phi)=\int_D
%\sigma(\mathbf{u}_s):\varepsilon(\mathbf{u}_s)\,\mathrm{d}x
%+\ell\int_D \chi(\phi)\,\mathrm{d}x\\
\begin{array}{cl}
& \begin{cases}
-\operatorname{div}\sigma(\mathbf{u}_s)=\mathbf{0} & \text{in}\ D,\\
\mathbf{u}_s=\mathbf{0} & \text{on}\ \Gamma_{\mathrm{fix}},\\
\sigma(\mathbf{u}_s)\,\mathbf{n}= {\mathbf{q}} & \text{on}\ \Gamma_{\mathrm{load}},\\
\sigma(\mathbf{u}_s)\,\mathbf{n}= \mathbf{0} & \text{on}\ \Gamma_{\mathrm{o}},
\end{cases}
\end{array}
\end{gathered}
\end{equation}
% \begin{equation}
% \begin{gathered}
% \min_{\phi}\ J(\phi)=\int_D
% \sigma(\mathbf{u}_s):\varepsilon(\mathbf{u}_s)\,\mathrm{d}x
% +\ell\int_D \chi(\phi)\,\mathrm{d}x\\
% \begin{array}{cl}
% \text{s.t.} & \begin{cases}
% -\operatorname{div}\sigma(\mathbf{u}_s)=\mathbf{0} & \text{in}\ D,\\
% \mathbf{u}_s=\mathbf{0} & \text{on}\ \Gamma_{\mathrm{fix}},\\
% \sigma(\mathbf{u}_s)\,\mathbf{n}= {\mathbf{q}} & \text{on}\ \Gamma_{\mathrm{load}},\\
% \sigma(\mathbf{u}_s)\,\mathbf{n}= \mathbf{0} & \text{on}\ \Gamma_{\mathrm{o}},
% \end{cases}
% \end{array}
% \end{gathered}
% \end{equation}
where $\mathbf{n}$ denotes the outward unit normal and the tangential
traction is given by $|{\mathbf{q}}|=10$. The isotropic constitutive
relation is
$$C\varepsilon(\mathbf{u}_s)=2\mu \varepsilon(\mathbf{u}_s) +\lambda{\rm div}\ \mathbf{u}_s I.$$
The tensor $I$ denotes the identity tensor, and the Lamé constants are
given by
\[
\lambda=\frac{E\nu}{(1+\nu)(1-2\nu)},\qquad
\mu=\frac{E}{2(1+\nu)},
\]
where $E=200$ is Young's modulus and $\nu=0.3$ is Poisson's ratio.
The selected parameter combination is reported in
Table~\ref{tab:params}; the corresponding displacement field is shown in
Figure~\ref{fig:representative}(b).

The initialization stage
took five rounds, i.e., one generation plus four repairs. The first repair round corrected two issues: a field-declaration error, in which the test function had been listed as a physical field with a type and unit, and a boundary-label error, in which the list of design-edge labels had been left empty; the second repair round corrected the declaration of the adjoint equation. The generation stage took three
rounds (one generation plus two repairs): the first repair round corrected
missing constant declarations, including the weak factor, the domain
sizes, and the Lamé constants; the second repair round corrected transcription
errors: the domain size and the weak factor did not match the
specification. Reflection passed the overall check and continued with
the post-processing.

The nonlinear variant replaces the material with a nonlinear
composite used in fibre-reinforced concrete: the stiffness
grows wherever the material deforms most, which helps the wall carry
the lateral load. Loads, supports, and the design setup are otherwise
identical
to the linear case; the only modification is the following
strain-dependent stiffness law
\begin{equation}
C(\varepsilon)=\bigl(1+\kappa\,\varepsilon:\varepsilon\bigr)C.
\end{equation}
Here, $\kappa=25$ is the strain-hardening coefficient.
As this strain-dependent law fails the major symmetry, no potential with
$\sigma$ as its gradient exists; the compliance term of the objective
therefore serves as a general performance measure. AutoFreeFEM's symbolic
differentiation generates the consistent tangent, including the derivative
of the stiffness factor with respect to $\varepsilon:\varepsilon$, and
assembles the adjoint from it. Table~\ref{tab:params} lists the parameter combination selected for this variant.

The pressure vessel end cap~\cite{szybinski2019flatends} is a flat vessel closure subjected to internal fluid
pressure, carrying the pressure to the welded
rim. As shown in Figure~\ref{fig:representative}(d), the design domain is a plane-strain strip between two symmetry
planes through the vessel axis. The pressurized top face
$\Gamma_{\mathrm{press}}$ and the fixed welded bottom rim
$\Gamma_{\mathrm{weld}}$ are prescribed; the two lateral symmetry
planes $\Gamma_{\mathrm{roll}}$, with their roller conditions, are
free for design. Consider \eqref{lsmTOP} with $\mathcal{J}(\phi)=\int_D
\sigma(\mathbf{u}_s):\varepsilon(\mathbf{u}_s)\,\mathrm{d}x$
subject to
\begin{equation}
\begin{array}{cl}
& \begin{cases}
-\operatorname{div}\sigma(\mathbf{u}_s)=\mathbf{0} & \text{in}\ D,\\
\mathbf{u}_s=\mathbf{0} & \text{on}\ \Gamma_{\mathrm{weld}},\\
\sigma(\mathbf{u}_s)\,\mathbf{n}=p\,\mathbf{n} & \text{on}\ \Gamma_{\mathrm{press}},\\
\sigma(\mathbf{u}_s)\,\mathbf{n}+k_p {u}_{sx}\mathbf{e}_x=\mathbf{0} & \text{on}\ \Gamma_{\mathrm{roll}}.
\end{cases}
\end{array}
\end{equation}
% \begin{equation}
% \begin{gathered}
% \min_{\phi}\ J(\phi)=\int_D
% \sigma(\mathbf{u}_s):\varepsilon(\mathbf{u}_s)\,\mathrm{d}x
% +\ell\int_D \chi(\phi)\,\mathrm{d}x\\
% \begin{array}{cl}
% \text{s.t.} & \begin{cases}
% -\operatorname{div}\sigma(\mathbf{u}_s)=\mathbf{0} & \text{in}\ D,\\
% \mathbf{u}_s=\mathbf{0} & \text{on}\ \Gamma_{\mathrm{weld}},\\
% \sigma(\mathbf{u}_s)\,\mathbf{n}=p\,\mathbf{n} & \text{on}\ \Gamma_{\mathrm{press}},\\
% \sigma(\mathbf{u}_s)\,\mathbf{n}+k_p {u}_{sx}\mathbf{e}_x=\mathbf{0} & \text{on}\ \Gamma_{\mathrm{roll}},
% \end{cases}
% \end{array}
% \end{gathered}
% \end{equation}
The ersatz constitutive law \eqref{fwijeo} is the same as in the
shear-wall case. Here, $E=200$ and $\nu=0.3$ are Young's modulus and
Poisson's ratio, $p=10$ is the pressure, $k_p=10^6$ is the roller
stiffness, $u_{sx}$ is the horizontal displacement component, and
$\mathbf{e}_x$ is the axial unit vector. The parameter combination selected for this case is listed in Table~\ref{tab:params}.

For the full run of this case, the initialization stage
took three rounds, i.e., one generation plus two repairs: the first repair round corrected a field-declaration problem and a unit-notation problem,
and the second repair round corrected the integrand of the volume term in the
objective. The generation stage took four rounds, i.e., one generation
plus three repairs: the first repair round corrected field-declaration,
objective-term, and load-sign problems, the second repair round corrected a
call-form problem and the Dirichlet notation, and the third repair round corrected
a missing Dirichlet label. Reflection fixed a few program issues. After
the re-execution, the overall check passed and the post-processing
proceeded.

The nonlinear variant uses the same end cap made of nonlinear
steel. Loads, supports, and the design setup are
otherwise identical to the linear case; the only modification is the
following strain-dependent stiffness law
\begin{equation}
C(\varepsilon)=\bigl(1+\kappa_h\,\varepsilon:\varepsilon\bigr)C.
\end{equation}
Here, $\kappa_h=50$ is the strain-hardening coefficient.
The compliance term of the objective serves as a general performance
measure, and the consistent tangent and adjoint are likewise produced by
AutoFreeFEM's symbolic differentiation. The corresponding parameter combination is reported in Table~\ref{tab:params}.

The four thermal diffusion cases, in contrast, are driven by uniform volumetric body heat. The cold plate is Bejan's classical volume-to-point benchmark~\cite{bejan1997constructal}: heat is
generated uniformly inside a power package, and all generated heat
must be conducted to a small cooling region where the heat sink is
attached. The design domain is
a 1\,m by 1\,m square plate with uniform volumetric heating and a
constant base conductivity.
The central 0.2\,m segment $\Gamma_{c}$
of the left edge, which is the cold patch, is maintained at the
reference temperature,
and every remaining boundary segment is insulated and free for design, as shown in Figure~\ref{fig:diffusion}(a).
The objective balances the dissipated thermal energy against the
material cost, and the optimization produces conductive paths from the
heated interior to the cold patch. More specifically, consider  \eqref{lsmTOP} with $\mathcal{J}(\phi)=\int_D
\chi(\phi)\,k\,|\nabla T|^2\,\mathrm{d}x$ subject to
\begin{equation}
\begin{gathered}
% \min_{\phi}\ J(\phi)=\int_D
% \chi(\phi)\,k\,|\nabla T|^2\,\mathrm{d}x+\ell\int_D \chi(\phi)\,\mathrm{d}x\\
\begin{array}{cl}
 & \begin{cases}
-\operatorname{div}\bigl(\chi(\phi)\,k\,\nabla T\bigr)=f & \text{in}\ D,\\
T=0 & \text{on}\ \Gamma_{c},\\
k\nabla T\cdot\mathbf{n}=0 & \text{on}\ \partial D\setminus\Gamma_{c},
\end{cases}
\end{array}
\end{gathered}
\end{equation}
with the volumetric heat source $f=10$ and conductivity $k=1$. The selected
parameter combination is reported in Table~\ref{tab:params}; the corresponding temperature field
is shown in Figure~\ref{fig:diffusion}(b).

For this case, the initialization stage
took four rounds (one generation plus three repairs): the three
repair rounds involved repeated attempts to correct errors in the
weak-material interpolation formula and the unit specifications. Reflection passed the overall check and continued with
the post-processing.

The nonlinear variant models a semiconductor plate whose conductivity decreases as the plate heats up, so that the
hottest regions are the least conductive. The geometry, the cold patch,
the uniform heat source, and the objective structure stay as in the
linear case; only the conductivity law changes, yielding a nonlinear thermal-diffusion model
\begin{equation}
k(T)=1/(1+0.1\,T).
\end{equation}
The selected parameter combination is listed in Table~\ref{tab:params}.

The two-sink cold plate represents a similar configuration~\cite{bejan1997constructal} with two cooling patches,
one on each side edge, instead of one: the heat generated throughout
the board must be conducted to the cooling patches, which are held at the
reference temperature while every remaining boundary segment is
insulated and free for design; the volumetric heat source and the
conductivity stay at $f=10$ and $k=1$. The selected parameter combination is reported in Table~\ref{tab:params}.

For this case, the initialization stage
took two rounds, i.e., one generation plus one repair. The repair round corrected
three issues: a mismatch between the field domain and the domain name, a
field-declaration error in which the test function had been listed as a
physical field with a type and unit, and a unit-notation error.
Reflection fixed a few program issues. After the re-execution, the
overall check passed and the post-processing proceeded.

\begin{figure}[!tb]
\centering
\includegraphics[width=0.98\textwidth]{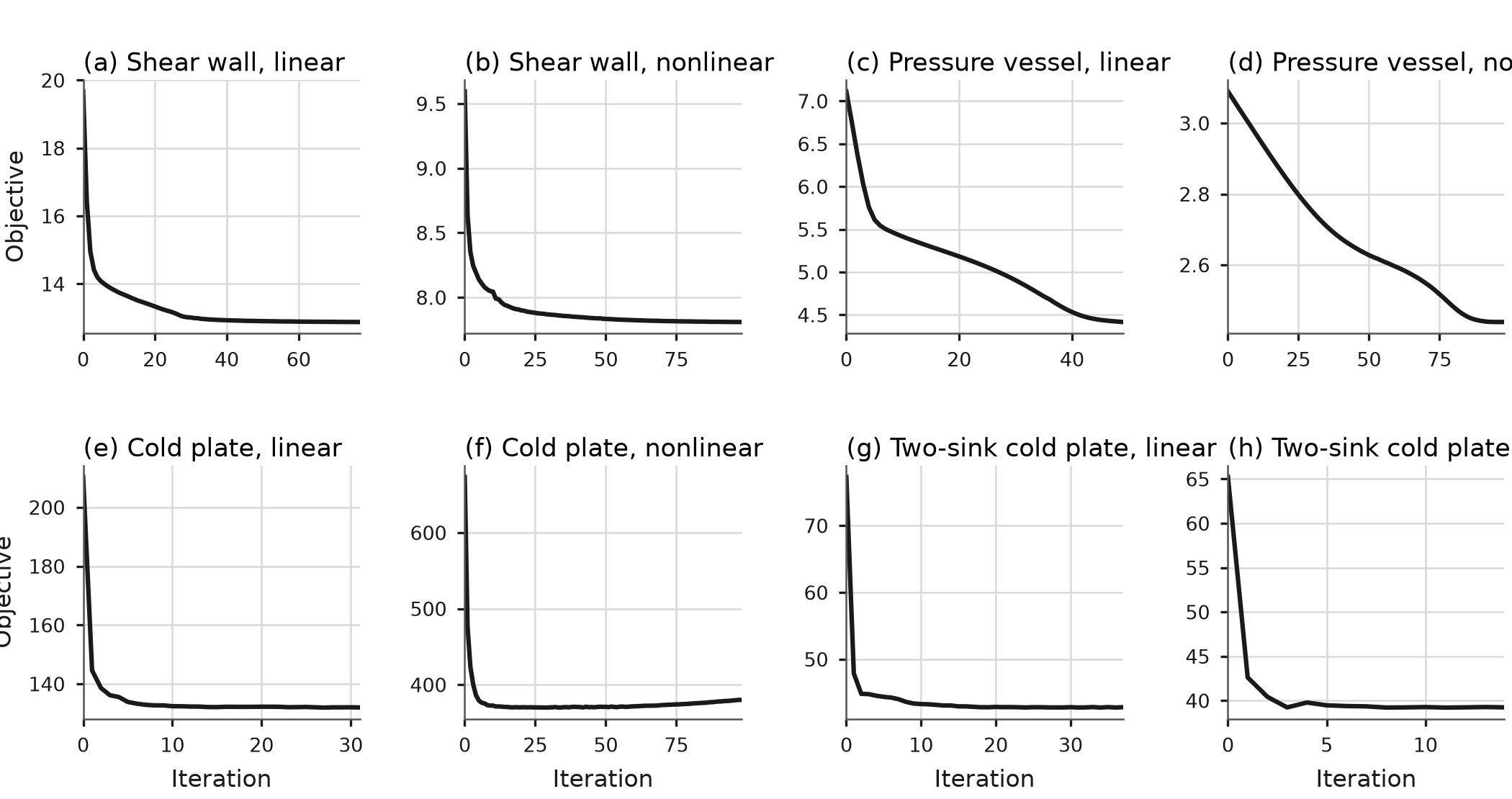}
\caption{Objective history of the eight cases:
(a) shear wall (linear); (b) shear wall (nonlinear); (c) pressure vessel
end cap (linear); (d) pressure vessel end cap (nonlinear); (e) cold
plate (linear); (f) cold plate (nonlinear); (g) two-sink cold plate (linear); (h) two-sink cold plate (nonlinear).}
\label{fig:objective}
\end{figure}

The nonlinear variant uses a polymer dielectric whose conductivity
increases with temperature. All other aspects of the two-sink
configuration remain unchanged; only the material law changes, as
follows
\begin{equation}
k(T)=1+0.1T.
\end{equation} 
The pair thus contrasts the nonlinear conduction law with the
linear one under the same geometry and heating. Table~\ref{tab:params} reports the parameter combination selected for this variant.
Figure~\ref{fig:objective} shows the objective history of the reported
experiment of each case. Within each history, the design attaining the
smallest objective value is the recorded design of that experiment.

\subsection{Ablation Study: the Role of the Pipeline Components}

We used a component ablation study to show how each component
contributed to the success rate and computational overhead. The initialization,
generation, and reflection stages share the same check-and-repair
logic, so we illustrate with the initialization stage here; the other
two stages lead to similar conclusions. For the eight representative cases of Section~4.2, in
addition to the normal pipeline, four variants were evaluated, each
missing one component: the stage prompts, the manual knowledge, the
reference examples, and
the repair history. Every experiment was subject to the same repair budget as the
normal version: one initial generation followed by at most ten
repairs, for a maximum of eleven rounds. The number of rounds and the
runtime were measured as described in Section~4.1; the normal-version data are the
initialization columns of Table~\ref{tab:runtime}.
Table~\ref{tab:initablation} collects the rounds, times, and outcomes
of the eight cases under the five configurations.

\begin{table}[!t]
\caption{Component ablation of the initialization stage on the eight
representative cases: the normal version, i.e., the initialization columns of
Table~\ref{tab:runtime}, against variants with one component removed.
Each cell is rounds/time; every experiment uses the same repair budget
as the normal version, i.e., one generation plus at most ten repairs.
$\dag$ marks a budget-exhausted failure; no\_examples,
no\_prompt,
no\_manual, and no\_history denote the variants missing the reference
examples, the stage prompts, the manual knowledge, and the repair
history, respectively.}
\label{tab:initablation}
\centering
\small
\begin{tabular}{l c c c c c}
\toprule
Case & Normal & no\_examples & no\_prompt & no\_manual & no\_history \\
\midrule
doubleSinkColdPlate & 2 / 46m32s & 3 / 20m02s & 11 / 35m27s\textsuperscript{$\dag$} & 11 / 60m51s\textsuperscript{$\dag$} & 8 / 46m53s \\
doubleSinkColdPlateNonlinear & 6 / 37m27s & 5 / 33m26s & 8 / 36m50s & 11 / 52m38s\textsuperscript{$\dag$} & 4 / 25m36s \\
pressureVessel & 3 / 20m32s & 4 / 29m22s & 11 / 76m55s\textsuperscript{$\dag$} & 11 / 44m13s\textsuperscript{$\dag$} & 6 / 37m22s \\
pressureVesselNonlinear & 2 / 11m00s & 11 / 68m21s\textsuperscript{$\dag$} & 10 / 30m53s & 11 / 42m40s\textsuperscript{$\dag$} & 11 / 67m05s\textsuperscript{$\dag$} \\
shearWall & 5 / 33m45s & 4 / 22m17s & 8 / 21m01s & 11 / 43m36s\textsuperscript{$\dag$} & 4 / 36m06s \\
shearWallNonlinear & 2 / 12m08s & 4 / 28m38s & 8 / 30m00s & 11 / 46m04s\textsuperscript{$\dag$} & 4 / 26m00s \\
vtpColdPlate & 4 / 30m54s & 4 / 25m11s & 10 / 35m48s & 8 / 36m56s & 3 / 22m05s \\
vtpColdPlateNonlinear & 7 / 68m18s & 11 / 85m14s\textsuperscript{$\dag$} & 8 / 29m34s & 10 / 35m26s & 11 / 74m11s\textsuperscript{$\dag$} \\
\midrule
Success rate & \textbf{8/8} & \textbf{6/8} & \textbf{6/8} & \textbf{2/8} & \textbf{6/8} \\
Average rounds & \textbf{3.9} & \textbf{5.8} & \textbf{9.3} & \textbf{10.5} & \textbf{6.4} \\
Average time & \textbf{32m35s} & \textbf{39m04s} & \textbf{37m03s} & \textbf{45m18s} & \textbf{41m55s} \\
\bottomrule
\end{tabular}
\end{table}

Table~\ref{tab:initablation} also shows the information each component
provides and the cost of losing it. The stage prompts provide a fixed
output format and structured repair procedures, thereby constraining
the open-ended task of structured problem modeling and simplifying the
reasoning required during
generation and repair; without them the outputs drift structurally.
The manual knowledge provides the domain-specific patterns of the specifications and the
physical criteria, thereby keeping the model outputs well-formed; without it,
the model struggles to produce a valid specification. The repair history provides
the record of the errors and the repair actions of the earlier rounds,
thereby keeping the repairer from repeating the same mistakes; without
it, the cases repeatedly oscillate between two or more error types and
struggle to reach a correct result. The reference examples provide
complete specifications of similar problems from the library, anchoring
the parser's design choices to proven patterns; without them the cases
took more rounds on average, and two nonlinear cases exhausted the
budget. Three concrete examples
illustrate the point. Without the manual knowledge, the pipeline failed to produce a valid
specification for the pressure vessel end cap because the $\phi$ field
remained empty across successive rounds; the case exhausted its 11-round budget and failed after 44m13s,
up from 20m32s. Without the reference examples, the nonlinear pressure
vessel end cap went from 2 rounds and 11m00s to 11 rounds and 68m21s,
exhausting the budget. Without the repair history, the two-sink cold plate
repeatedly made the same class of mistakes, taking eight rounds
instead of two rounds before it finally passed. The removal of repair history was associated with more repair rounds,
longer runtimes, and a lower success rate across the evaluated cases.

For the same task, the full pipeline achieved a
100\% success rate with the lowest average number of rounds and the
shortest average runtime;
removing any component always reduced the success rate and/or increased
the number of invalid attempts and the repair overhead. These results demonstrate the value of the pipeline components in
improving reliability and reducing repair overhead by providing
structured guidance, domain knowledge, and historical feedback.

\subsection{Error-Injection Experiments}

We used error-injection experiments to evaluate the pipeline's
ability to repair each class of error. As in the ablation study, we
focus on the generation stage; the other two stages lead to similar
conclusions. For the
eight representative cases of Section~4.2, we
injected one representative defect from each of the six injectable
error classes of Table~\ref{tab:genchecks} (program structure, call
forms, transcription, weak-form expressibility, penalty rules, and
mathematical equivalence) into a copy of the case code, then started the
generation repair flow, which continued until the generated code
passed the generation-stage validation checks or the repair budget
was exhausted. Table~\ref{tab:errorinjection}
collects the repair rounds and times of the 48 experiments.

\begin{table}[!t]
\caption{Error-injection experiments: one representative defect from
each of the six injectable error classes of Table~\ref{tab:genchecks} is
manually injected into a copy of the generated code of the eight
representative cases, and the generation repair flow is started with
the correct inputs. Each cell is repair rounds/time; the repair budget
is at most ten repairs, and all 48 experiments pass.}
\label{tab:errorinjection}
\centering
\footnotesize
\setlength{\tabcolsep}{4pt}
\begin{tabular}{l c c c c c c}
\toprule
Case & Structure & Call form & Transcription & Weak form & Penalty & Math equivalence \\
\midrule
doubleSinkColdPlate & 2 / 18m29s & 1 / 11m14s & 1 / 9m02s & 1 / 8m21s & 1 / 9m37s & 1 / 6m04s \\
doubleSinkColdPlateNonlinear & 1 / 10m42s & 6 / 52m18s & 1 / 13m19s & 1 / 13m13s & 2 / 19m47s & 2 / 23m00s \\
pressureVessel & 1 / 12m11s & 2 / 40m12s & 1 / 13m44s & 1 / 13m21s & 2 / 19m09s & 1 / 20m47s \\
pressureVesselNonlinear & 2 / 27m54s & 2 / 28m42s & 2 / 28m06s & 2 / 24m22s & 2 / 26m48s & 2 / 26m20s \\
shearWall & 1 / 18m21s & 2 / 18m59s & 1 / 13m04s & 2 / 17m42s & 1 / 13m56s & 1 / 10m50s \\
shearWallNonlinear & 2 / 16m08s & 1 / 8m42s & 1 / 11m47s & 2 / 17m24s & 1 / 10m16s & 1 / 7m43s \\
vtpColdPlate & 1 / 5m20s & 1 / 7m00s & 1 / 6m15s & 1 / 8m41s & 1 / 6m26s & 1 / 6m43s \\
vtpColdPlateNonlinear & 3 / 17m54s & 1 / 11m59s & 1 / 10m30s & 2 / 18m31s & 1 / 9m45s & 1 / 7m25s \\
\midrule
Success rate & \multicolumn{6}{c}{\textbf{100\%}} \\
Average rounds & \textbf{1.6} & \textbf{2.0} & \textbf{1.1} & \textbf{1.5} & \textbf{1.4} & \textbf{1.2} \\
Average time & \textbf{15m52s} & \textbf{22m23s} & \textbf{13m13s} & \textbf{15m12s} & \textbf{14m28s} & \textbf{13m36s} \\
\bottomrule
\end{tabular}
\end{table}

Table~\ref{tab:errorinjection} shows that all 48 injected-error
experiments passed the generation-stage validation checks within the
repair budget, yielding a 100\% success rate. The deterministic
checks caught the program-structure, call-form, transcription,
weak-form expressibility, and penalty-rule errors directly, while the
semantic review caught the mathematical-equivalence errors, so
the two verification layers provided complementary error detection capabilities; every error class was repaired in
1.1--2.0 rounds and 13m13s--22m23s on average, with call-form errors
the slowest, but still far below the budget. The pipeline thus repaired
every class of generation error quickly.

\subsection{Extension Cases: Fluid Flow and Thermo-Elastic Coupling}

We also tested other types of optimization problems.
The two extension cases
below include a fluid-flow problem and a coupled thermo-elastic problem, and
both were generated solely from natural-language
descriptions.

The fluid case is the level-set recast of the double-pipe benchmark~\cite{farrell2013dolfinadjoint}, a topology
optimization problem in Stokes flow~\cite{borrvall2003topology}: the
objective is to design channels connecting the two
inlet regions on the left boundary to the two outlet regions on the
right boundary while balancing energy dissipation against the
fluid-volume
penalty. The
design domain is a 1.5 by 1 rectangle filled with a fluid of unit
viscosity $\mu=1$. The
solid acts on the flow through a Darcy resistance whose inverse
permeability is given by $\alpha(\chi)=\alpha_s+(\alpha_l-\alpha_s)\chi$, with the fluid
value $\alpha_l=2.5\times10^{-4}$ and the solid value
$\alpha_s=2.5\times10^{4}$. On the left edge, two strips of width $1/6$,
centered at $y=0.25$ and $y=0.75$, are assigned parabolic velocity
profiles with a peak value of 1 at the center of each strip; the two
outlet strips on the right edge are left free for the outflow, and the
remaining boundary is treated as a no-slip wall. The initial mesh
contained 3{,}662 nodes.

\begin{figure}[!tb]
\centering
\subfloat{%
  \begin{minipage}{0.54\textwidth}
  \centering
  \includegraphics[width=0.9\linewidth]{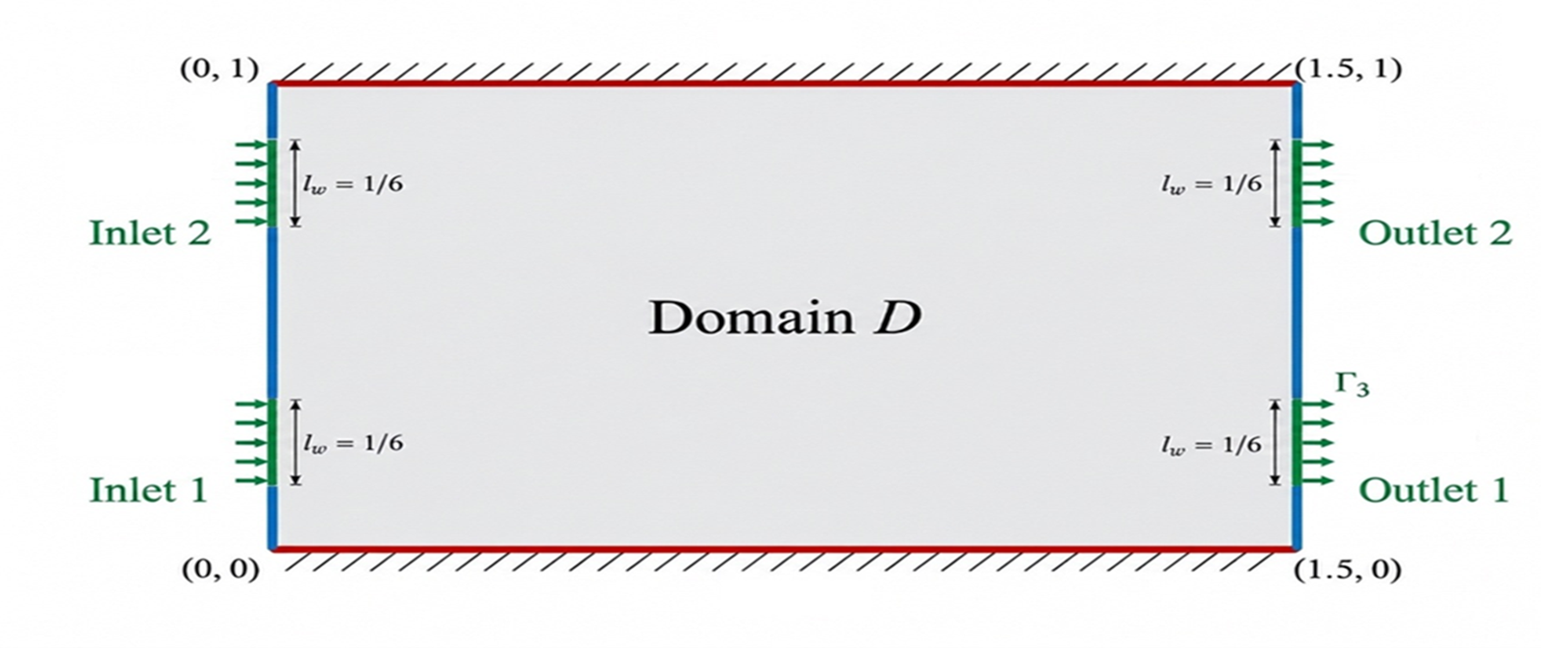}
  \end{minipage}}\\[0.5ex]
\subfloat{%
  \begin{minipage}{0.49\textwidth}
  \centering
  \includegraphics[width=0.8\linewidth]{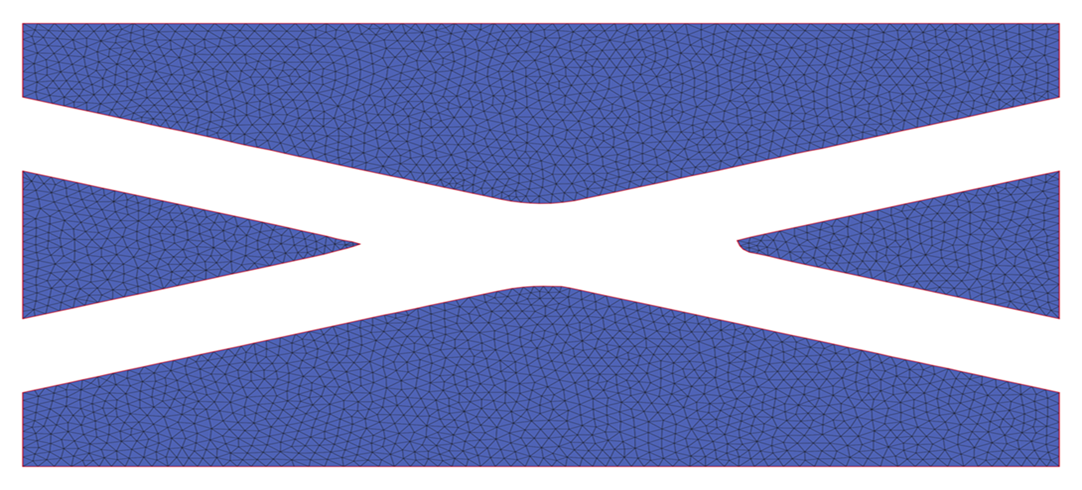}
  \end{minipage}}\hfill
\subfloat{%
  \begin{minipage}{0.49\textwidth}
  \centering
  \includegraphics[width=0.8\linewidth]{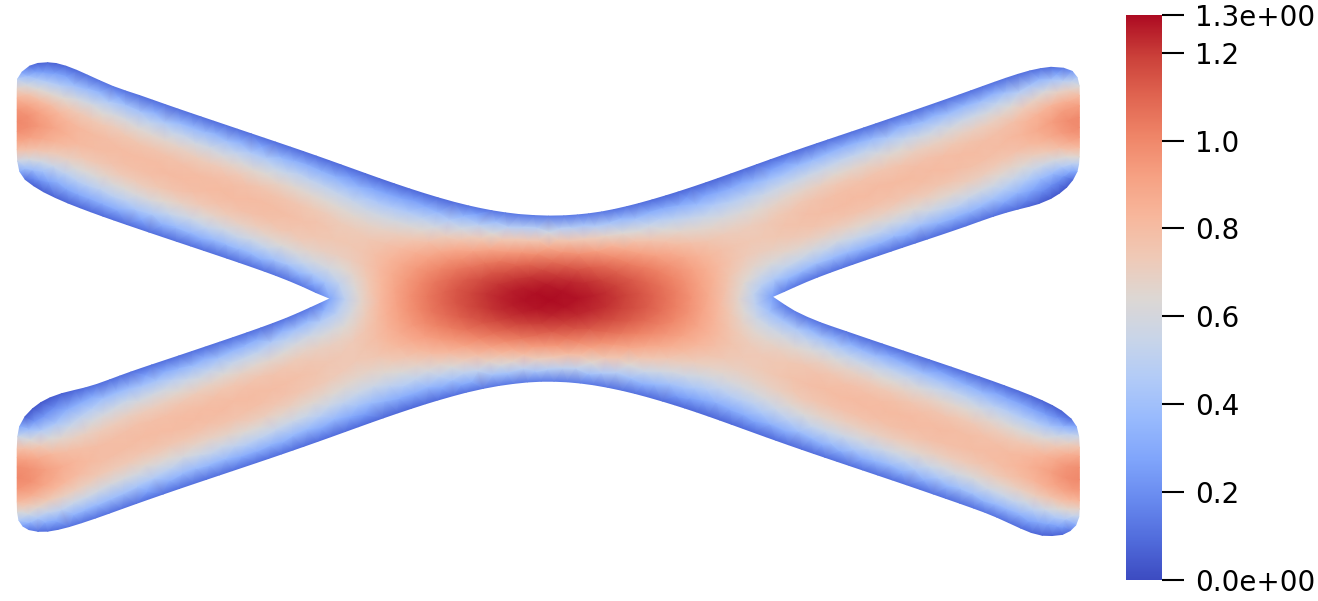}
  \end{minipage}}
\caption{Topology optimization of Stokes flow: problem
schematic (top); initial design with the finite-element mesh (bottom left); and the velocity field of the optimized design (bottom right).}
\label{fig:stokestopology}
\end{figure}

The objective combines the dissipated power, consisting of Darcy
friction and viscous dissipation, with a material-cost term that
controls the fluid volume: here the level-set function designates the
fluid channel, so that $\chi=1$ in the fluid region and the volume
functional $\mathcal{V}$ measures the fluid volume. Consider \eqref{lsmTOP} with~$\mathcal{J}(\phi)=\frac{1}{2}\int_D (\alpha(\chi)\,\lvert\mathbf{u}\rvert^2
+\mu|\nabla\mathbf{u}|^2)\,\mathrm{d}x$~\mbox{subject to}
\begin{equation}
\begin{gathered}
\begin{array}{cl}
& \begin{cases}
-\mu\Delta\mathbf{u}+\alpha(\chi)\,\mathbf{u}+\nabla p=\mathbf{0} & \text{in}\ D,\\
\operatorname{div}\mathbf{u}=0 & \text{in}\ D,\\
\mathbf{u}=\mathbf{b} & \text{on}\ \Gamma_{\mathrm{in}},\\
\sigma(\mathbf{u},p)\,\mathbf{n}=\mathbf{0} & \text{on}\ \Gamma_{\mathrm{out}},\\
\mathbf{u}=\mathbf{0} & \text{on}\ \Gamma_{\mathrm{w}},
\end{cases}
\end{array}
\end{gathered}
\end{equation}
where $\mathbf{u}$ is the flow velocity, $p$ is the pressure,
$\Gamma_{\mathrm{in}}$ and $\Gamma_{\mathrm{out}}$ denote the inlet
strips on the left edge and the outlet strips on the right edge,
respectively, $\mathbf{b}$ denotes the prescribed parabolic strip
profile, $\sigma(\mathbf{u},p)=\mu(\nabla\mathbf{u}+\nabla\mathbf{u}^{T})-pI$ is the stress tensor, $\ell=0.05$ is the material-cost weight, and
$\Gamma_{\mathrm{w}}=\partial D\setminus(\Gamma_{\mathrm{in}}\cup\Gamma_{\mathrm{out}})$
is the remaining no-slip wall. The optimization produced
the bifurcated double-pipe design of Figure~\ref{fig:stokestopology}: two
arms from the inlet strips merge in a wide central chamber and fork
again to the outlet strips.

The thermo-elastic case~\cite{xia2008thermoelastic} models the
plane-strain cross-section of a long structure defined on a
4\,m by 1\,m rectangular domain. The solid has a Young's modulus of
32\,MPa, a Poisson's ratio of 0.23, a density of
$\rho=2500$\,kg/m$^3$, a coefficient of linear thermal expansion of
$\alpha_T=3.6\times10^{-5}\,\mathrm{K}^{-1}$, and a thermal conductivity
of $k_0=1.25$\,W/(m\,K). On the left
edge $\Gamma_{T}$, a temperature profile is prescribed and imposed
weakly through a Robin penalty term. The prescribed temperature
increases linearly from 0\,$^\circ$C at the bottom corner to
$30\,^\circ$C at the top corner. The left and
right edges $\Gamma_{D}$ are clamped by penalty terms; the top edge
$\Gamma_{N}$ carries a uniform downward traction of magnitude
$q=5000\,\mathrm{N/m^2}$;
gravity $g=9.81$\,m/s$^2$ acts downward on every solid region; and the
bottom edge $\Gamma_b$ is insulated and traction-free. The initial mesh
contained 10{,}251 nodes.

\begin{figure}[!tb]
\centering
\subfloat{%
  \begin{minipage}{0.49\textwidth}
  \centering
  \includegraphics[width=\linewidth]{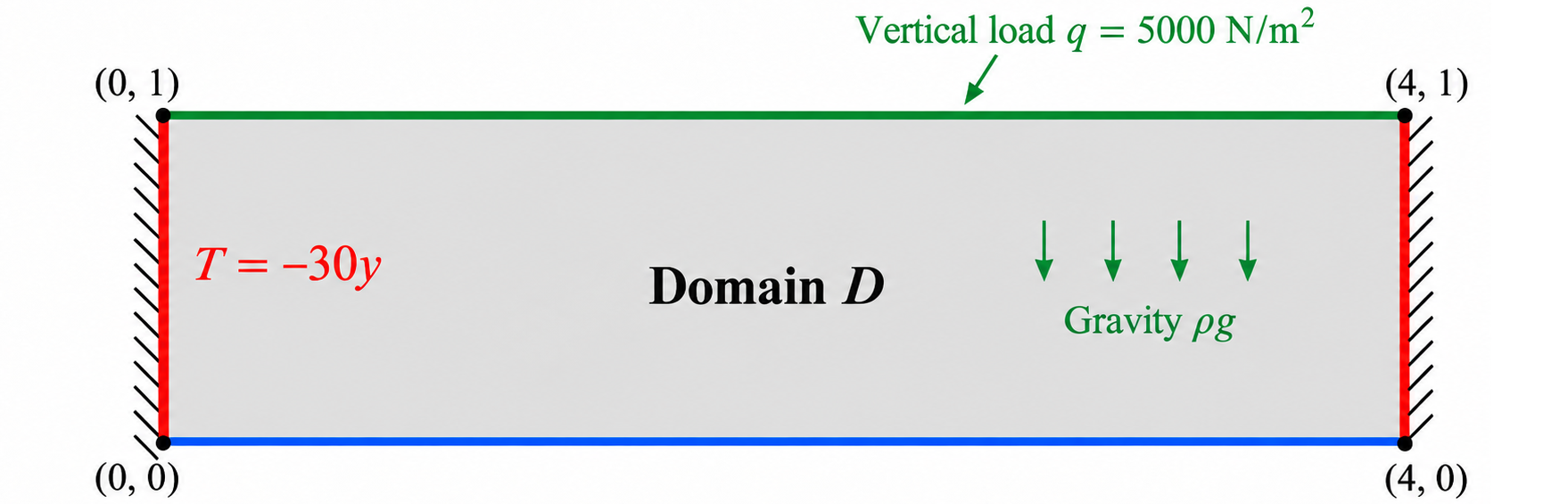}
  \end{minipage}}\hfill
\subfloat{%
  \begin{minipage}{0.49\textwidth}
  \centering
  \includegraphics[width=\linewidth]{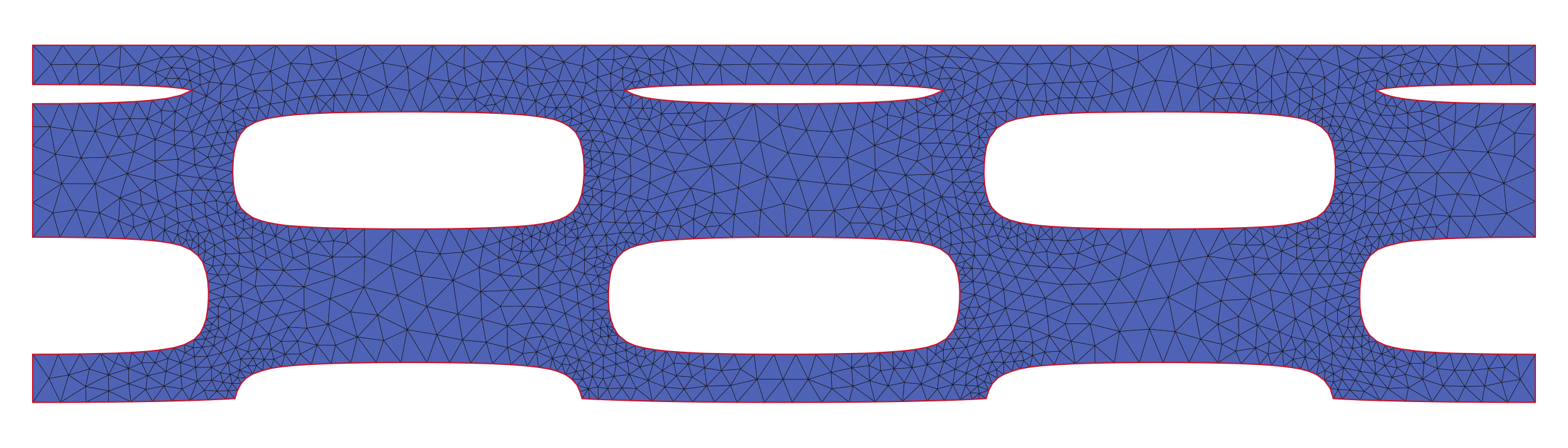}
  \end{minipage}}\\[0.5ex]
\subfloat{%
  \begin{minipage}{0.49\textwidth}
  \centering
  \includegraphics[width=\linewidth]{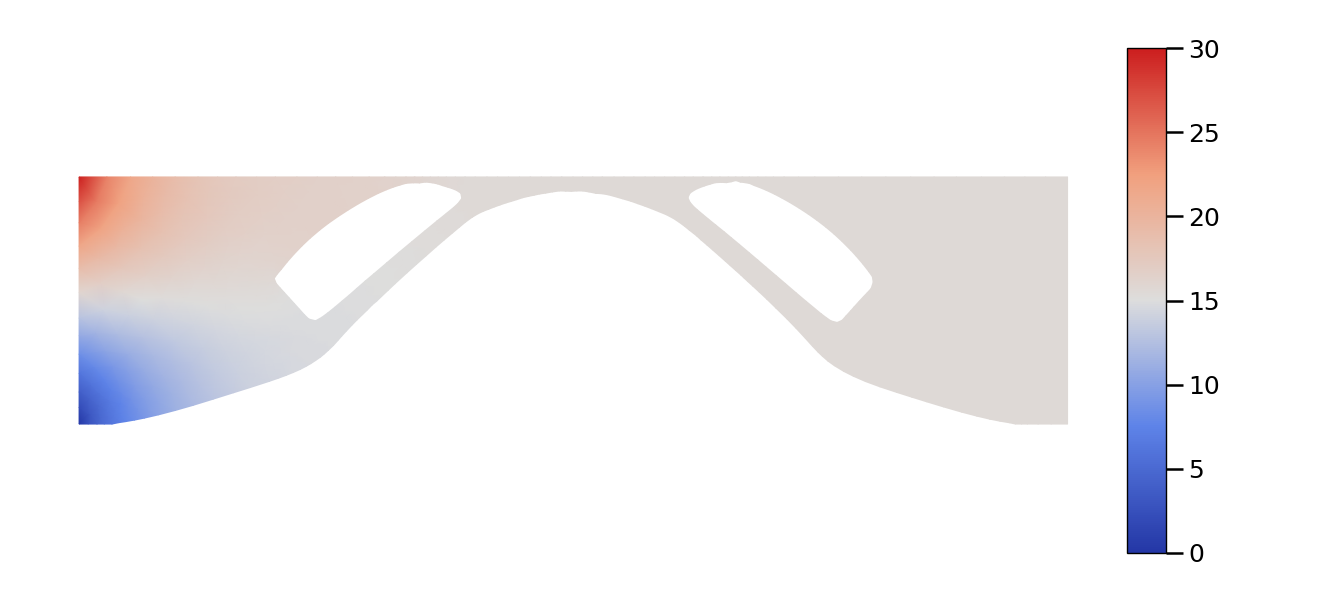}
  \end{minipage}}\hfill
\subfloat{%
  \begin{minipage}{0.49\textwidth}
  \centering
  \includegraphics[width=\linewidth]{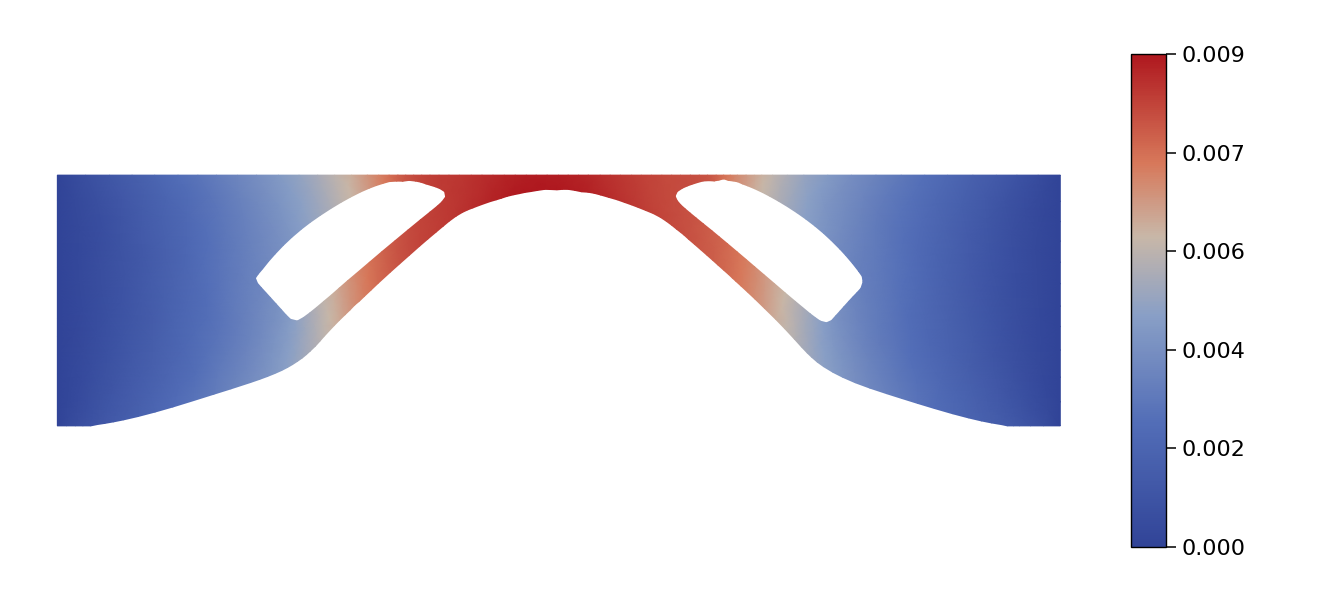}
  \end{minipage}}
\caption{The thermo-elastic case: problem schematic (top left);
initial design (top right); the optimized design colored by the
temperature field $T$ (bottom left); and the optimized design colored
by the displacement field $\mathbf{u}$ (bottom right).}
\label{fig:thermo-elastic}
\end{figure}

The objective combines the material cost, a penalty on the
horizontal displacement of the clamped edges, and a term that
minimizes the downward displacement of the loaded edge. More
specifically, consider \eqref{lsmTOP} with~$\mathcal{J}(\phi)=\beta\int_{\Gamma_{D}} \left(\mathbf{e}_x\cdot\mathbf{u}\right)^2\,\mathrm{d}s
-\gamma\int_{\Gamma_{N}} \left(\mathbf{e}_y\cdot\mathbf{u}\right)\,\mathrm{d}s$ subject to
\begin{equation}
\begin{gathered}
% \begin{aligned}
% \min_{\phi}\ J(\phi)=\beta&\int_{\Gamma_{D}} \left(\mathbf{e}_x\cdot\mathbf{u}\right)^2\,\mathrm{d}s\\
% -\gamma&\int_{\Gamma_{N}} \left(\mathbf{e}_y\cdot\mathbf{u}\right)\,\mathrm{d}s+\ell\int_D \chi(\phi)\,\mathrm{d}x
% \end{aligned}\\
\begin{array}{cl}
& \begin{cases}
-\operatorname{div}\sigma=-\chi(\phi)\,\rho g\,\mathbf{e}_y & \text{in}\ D,\\
\sigma=C:\bigl(\varepsilon(\mathbf{u})-\alpha_T T I\bigr) & \text{in}\ D,\\
-\operatorname{div}\bigl(k\,\nabla T\bigr)=0 & \text{in}\ D,\\
k\nabla T\cdot\mathbf{n}+\beta_T(T-30y)=0 & \text{on}\ \Gamma_{T},\\
\sigma\,\mathbf{n}=-q\,\mathbf{e}_y & \text{on}\ \Gamma_{N},\\
\mathbf{u}=\mathbf{0} & \text{on}\ \Gamma_{D},\\
k\nabla T\cdot\mathbf{n}=0 & \text{on}\ \partial D\setminus\Gamma_{T},\\
\sigma\,\mathbf{n}=\mathbf{0} & \text{on}\ \Gamma_{b},
\end{cases}
\end{array}
\end{gathered}
\end{equation}
with $\beta=10^{11}$, $\beta_T=10^{11}$, and $\gamma=100$. The stiffness $C$, the thermal expansion coefficient $\alpha_T$, and the conductivity $k$ are material-interpolated: each equals its solid value given above in the solid region and one thousandth of that value in the void. Here, $\mathbf{u}=\mathbf{0}$ on $\Gamma_{D}$ is an ideal boundary condition, imposed approximately through a penalty term with coefficient $\beta$ rather than enforced strictly; the same coefficient $\beta$ also weights the displacement-penalty term in the objective. We set $\ell=1$.
The optimization converged to the design of
Figure~\ref{fig:thermo-elastic}: the material concentrates in a
load-carrying band that spans the domain between the clamped edges
and is perforated where the load path is weak; the temperature field
varies across the structure from the heated left edge.

% =====================================================================
\section{Conclusion}

We have presented Agent4FreeFem, a pipeline that turns a
natural-language description of a topology optimization problem into a
verified, optimized design without expert intervention. The pipeline
chains four stages, initialization, generation, execution, and
reflection, under one methodological thread: the problem definition
stays unmodified, errors are localized using independent sources of
evidence, and constrained repairs are completed automatically, making
the pipeline efficient, end-to-end, autonomous, and verifiable. The
evaluation encompasses both single-physics and coupled multiphysics
problems, and every case succeeded. The results demonstrate that the
pipeline is a practical step toward making topology optimization
more accessible by reducing the expertise required for solver-code development.

% =====================================================================
\section*{Acknowledgments}
The authors would like to thank Olivier Pironneau, Fr\'ed\'eric Hecht,
and Pierre-Henri Tournier of the Laboratoire Jacques-Louis Lions,
Sorbonne Universit\'e, and Ahmed Chabib of Institut Polytechnique de
Paris for their inspiration and guidance in this work.

This work was supported in part by the National Natural Science Foundation of China under grant
Nos. 12471377 and U25A20200, and by the Science and Technology Commission of Shanghai
Municipality under grant No. 22DZ2229014.

\section*{Data Availability}
The code supporting this study is available at
\url{https://github.com/PierreCUI/Agent4Optimization}.

\bibliographystyle{unsrt}
\bibliography{references}

\end{document}